\documentclass[aap,MSNbibl,seceqn,dvips]{arximspdf}
\usepackage{mathrsfs}
\usepackage{graphicx}

\doi{10.1214/12-AAP888} 
\volume{23}
\issue{5}
\pubyear{2013}
\firstpage{1841}
\lastpage{1878}

\makeatletter
\newcommand{\varliminf}{\mathop{\underline{\lim}}}
\newcommand{\varlimsup}{\mathop{\overline{\lim}}}
\newtheorem{thmm}{Theorem}[section]
\newtheorem{prop}{Proposition}[section]
\newtheorem{cor}{Corollary}[section]
\newtheorem{lem}{Lemma}[section]

\newcommand{\aaa}{C_{4}}
\newcommand{\bbb}{C_{6}}
\newcommand{\ccc}{C_{12}}
\newcommand{\ddd}{C_{14}}
\newcommand{\eee}{C_{1}}
\newcommand{\fff}{C_{2}}
\newcommand{\hhh}{C_{3}}
\newcommand{\iii}{C_{7}}
\newcommand{\jjj}{C_{8}}
\newcommand{\kkk}{C_{9}}
\newcommand{\mmm}{C_{11}}
\newcommand{\nnn}{C_{10}}
\newcommand{\ooo}{C_{13}}
\newcommand{\doo}{C_{5}}

\newcommand{\aba}{a_{1}}
\newcommand{\aca}{a_{2}}

\makeatother

\begin{document}
\begin{frontmatter}

\title{Decay of tails at equilibrium for FIFO join the shortest queue networks}
\runtitle{Join the shortest queue}

\begin{aug}
\author[A]{\fnms{Maury}~\snm{Bramson}\corref{}\thanksref{t1}\ead[label=e1]{bramson@math.umn.edu}},
\author[B]{\fnms{Yi}~\snm{Lu}\ead[label=e2]{yilu4@illinois.edu}}
\and
\author[C]{\fnms{Balaji}~\snm{Prabhakar}\thanksref{t3}\ead[label=e3]{balaji@stanford.edu}}
\thankstext{t1}{Supported in part by NSF Grants CCF-0729537 and DMS-11-05668.}
\thankstext{t3}{Supported in part by NSF Grant CCF-0729537 and by a
grant from the Clean Slate Program at Stanford University.}
\runauthor{M. Bramson, Y. Lu and B. Prabhakar}
\affiliation{University of Minnesota, University of Illinois and
Stanford University}
\address[A]{M. Bramson\\
University of Minnesota\\
Twin Cities Campus\\
School of Mathematics\\
206 Church Street S.E.\\
Minneapolis, Minnesota 55455\\
USA\\
\printead{e1}}

\address[B]{Y. Lu\\
University of Illinois\\
Coordinated Science Lab\\
1308 W. Main Street\\
Urbana, Illinois 61801\\
USA\\
\printead{e2}}

\address[C]{B. Prabhakar\\
Stanford University\\
David Packard Building, Room 269\\
Stanford, California 94305\\
USA\\
\printead{e3}}
\end{aug}

\received{\smonth{6} \syear{2011}}
\revised{\smonth{4} \syear{2012}}

%
\begin{abstract}
In join the shortest queue networks, incoming jobs
are assigned to the shortest queue from among a randomly chosen subset
of $D$ queues, in a system of $N$ queues;
after completion of service at its queue, a
job leaves the network. We also assume that jobs arrive into the system
according to a
rate-$\alpha N$ Poisson process, $\alpha< 1$, with rate-$1$ service at
each queue.
When the service at queues is exponentially distributed, it was shown in
Vvedenskaya et al.
[\textit{Probl. Inf. Transm.} \textbf{32} (1996) 15--29]
that the tail of the equilibrium
queue size decays doubly exponentially in the limit as
$N\rightarrow\infty$. This is a substantial improvement over the case
$D=1$, where the queue size decays exponentially.\looseness=-1

The reasoning in [\textit{Probl. Inf. Transm.} \textbf{32} (1996) 15--29]
does not easily generalize to jobs with nonexponential
service time distributions.
A modularized program for treating general service time distributions
was introduced
in Bramson et al. [In \textit{Proc. ACM SIGMETRICS} (2010)
275--286]. The program relies on an ansatz that asserts, in
equilibrium, any fixed number of queues become independent of one another
as $N\rightarrow\infty$. This ansatz was demonstrated in several
settings in Bramson
et al.
[\textit{Queueing Syst.} \textbf{71} (2012) 247--292],
including for networks where the service discipline is FIFO and the service
time distribution has a decreasing hazard rate.

In this article, we investigate the limiting behavior, as $N\rightarrow
\infty$,
of the equilibrium at a queue when
the service discipline is FIFO and the service time distribution has a
power law with
a given exponent $-\beta$, for $\beta> 1$.
We show under the above ansatz that, as $N\rightarrow\infty$, the
tail of the equilibrium queue size
exhibits a wide range of behavior depending on the relationship between
$\beta$ and $D$.
In particular, if $\beta> D/(D-1)$, the tail is doubly exponential
and, if
$\beta< D/(D-1)$, the tail has a power law. When $\beta= D/(D-1)$,
the tail is
exponentially distributed.
\end{abstract}

%

\end{frontmatter}

\section{Introduction}\label{sec1}
We consider \emph{join the shortest queue} (JSQ) networks, where
incoming ``jobs'' (or ``customers'')
are assigned to the shortest queue from among~$D$ distinct queues, $D
\ge2$, with these queues being chosen uniformly
from among the $N$ queues in the system, with $D\le N$. When two or
more of these queues each have the
fewest number of jobs, each of the queues is chosen with equal probability.
After completion of service at its queue, a
job leaves the network. We assume that jobs arrive according to a
rate-$\alpha N$ Poisson process, $\alpha< 1$, and that jobs are served
independently and at rate $1$ at each queue.
We are interested in this article in the case where the service discipline
at each queue is first-in, first-out (FIFO).

When the service at queues is exponentially distributed, the evolution
of the system is given by
a countable state Markov chain where a state is given by the number of
jobs at each queue. It
is not difficult to show that a unique equilibrium distribution exists;
this equilibrium is exchangeable with respect to the ordering of the queues.
Let $P_k^{(N)}$ denote the probability that there are at least $k$ jobs
in equilibrium
for the system with $N$ queues.
It was shown in Vvedenskaya et al.~\cite{VDK} that
%
%
\begin{equation}
\label{eq111} \lim_{N\rightarrow\infty}P_k^{(N)} =
\alpha^{(D^{k} -1)/(D-1)}\qquad \mbox{for } k \in\mathbb{Z}_+;
\end{equation}
%
%
in particular, the right tail of $P_k^{(N)}$ decays doubly exponentially
fast in the limit as $N\rightarrow\infty$.
This behavior is a substantial improvement over the case
$D=1$, where $P_k^{(N)}$ decays exponentially, and has led to
substantial interest in JSQ networks in the literature. For other
references, see
Azar et~al.~\cite{ABKU}, Graham~\cite{G}, Luczak--McDiarmid \cite
{LM1,LM2}, Martin--Suhov~\cite{MS},
Mitzenma-\break cher~\cite{M}, Suhov--Vvedenskaya~\cite{SV}, Vocking~\cite{V}
and Vvedenskaya--Suhov~\cite{VS}.

Little work has been done on the behavior of JSQ networks when the
service times are not
exponentially distributed. In this setting, the underlying Markov
process will typically have
an uncountable state space, and positive Harris recurrence for the
process is no longer obvious. The latter
was shown in Foss--Chernova~\cite{FC}, and uniform bounds on the
equilibria were shown in
Bramson~\cite{B2}. (Both articles also considered JSQ networks with
more general arrivals and routing
of jobs.)

This paper builds on previous work~\cite{B2,BLP1} and \cite
{BLP2} by the authors.
Bramson et al.~\cite{BLP1} described a modularized program for
analyzing the limiting behavior
of the equilibria of a family of JSQ networks with general service
times, as $N\rightarrow\infty$.
An important step is to show that any fixed number of queues become
independent of one
another, with each converging to a limiting distribution that is
the equilibrium for an associated Markov process with a single queue,
which is a \emph{cavity process}. This process corresponds, in an
appropriate sense, to
``setting $N=\infty$'' in the JSQ network and
viewing the corresponding infinite dimensional process at a single
queue. We will
refer to this equilibrium as the \emph{equilibrium environment}. In
Section~\ref{sec2}, we will precisely
define this terminology.

Although it seems that this independence should hold in a very general
setting, including
under a wide range of service disciplines, demonstrating it appears to
be a difficult problem. In
Bramson et al.~\cite{BLP1}, this independence and convergence to the
equilibrium environment
were stated as an ansatz. This ansatz was demonstrated in Bramson et
al.~\cite{BLP2}
in several settings including for networks where the service discipline
is FIFO and the service
distribution has a decreasing hazard rate.

In this article, we employ the restriction of the above ansatz to FIFO
networks. This version of the ansatz will be
precisely stated in Section~\ref{sec2}. Here, we summarize it for application in
the current section:
%
%
 %
\begin{equation}\label{eq112}
\begin{tabular}{p{300pt}@{}}
For a family of networks with the FIFO service
discipline that are all in equilibrium, any fixed number of queues become
independent in the limit as $N\rightarrow\infty$. Moreover, each
marginal
distribution converges to the unique associated equilibrium environment.
%
\end{tabular}
\end{equation}
Although this ansatz has only been
demonstrated for service distributions having decreasing hazard rate
and for general service distributions when the arrival rate $\alpha$
is sufficiently small,
our arguments here do not otherwise require either restriction. Other
applications
of the ansatz, but for the processor sharing and LIFO service
disciplines, are given in~\cite{BLP1}.

Our goal, in this article, is to investigate the limiting behavior
of the right tail of the associated equilibrium environment,
under the FIFO service discipline and with the assigned mean-$1$
service distribution
$F(\cdot)$. Denote by $P_k$ the probability that there are at least
$k$ jobs in the equilibrium
environment. We will show that, when $F(\cdot)$ has a power law tail with
exponent $-\beta$, for given $\beta> 1$, the tail of $P_k$
exhibits a wide range of behavior depending on the relationship between
$\beta$ and $D$.
In particular, if $\beta> D/(D-1)$, the tail is doubly exponential
and, if
$\beta< D/(D-1)$, the tail has a power law; when $\beta= D/(D-1)$,
the tail is
exponentially distributed. When $\beta\nearrow\infty$, the
coefficient $q_D(\beta)$ of $k$
in the doubly exponential tail
converges to $1$, which is the coefficient of $k$ in (\ref{eq111}).
One obtains the same coefficient of $k$ whether $F(\cdot)$ has an
exponential tail or
has bounded support. Our main results are Theorems~\ref{thm121}, \ref
{thm131} and
\ref{thm141}. Theorem~\ref{thm121} covers the case $\beta> D/(D-1)$,
Theorem~\ref{thm131} covers
the case $\beta< D/(D-1)$ and Theorem~\ref{thm141} covers the case
$\beta= D/(D-1)$.
We set $\bar{F}(s) = 1 - F(s)$.

%
\begin{thmm}
\label{thm121}
Consider a family of JSQ networks, with given $D\ge2$ and
$N=D,D+1,\ldots,$ where the $N$th network
has Poisson rate-$\alpha N$ input, with $\alpha< 1$, and where service
at each queue is FIFO, with
distribution $F(\cdot)$ having mean $1$. Assume that (\ref{eq112})
holds and that
%
%
\begin{equation}
\label{eq122} %
\lim_{s\rightarrow\infty} \log{\bar{F}(s)}/\log{s} = -\beta,\vadjust{\goodbreak}
\end{equation}
with $\beta\in(D/(D-1),\infty)$. Then,
%
%
\begin{equation}
\label{eq123} %
\lim_{k\rightarrow\infty}(1/k)\log_D{
\log{(1/P_k)}} = q_D(\beta) %
\end{equation}
for some $q_D(\beta) \in(0,1)$. Moreover, $q_D(\beta)$ is continuous
in $\beta$ and
%
%
\begin{equation}
\label{eq124} %
q_D(\beta)\nearrow1\qquad  \mbox{exponentially
fast as } \beta\nearrow\infty. 
%
\end{equation}
When (\ref{eq122}) holds with $\beta= \infty$, then
(\ref{eq123}) holds with $q_D(\infty) = 1$.
\end{thmm}
Theorem~\ref{thm121} implies that,
when $\bar{F}(s) \sim cs^{-\beta}$ as $s\rightarrow\infty$,
for $\beta\in(D/\break(D-1),\infty)$ and $c>0$,
then $P_k = \exp\{-D^{(1+o(1))q_D(\beta)k}\}$.

%
\begin{thmm}
\label{thm131}
Consider a family of JSQ networks as in Theorem~\ref{thm121}, with~(\ref{eq122}) instead holding for
$\beta\in(1,D/(D-1))$. Then
%
%
\begin{equation}
\label{eq132} %
\lim_{k\rightarrow\infty} \log{(1/P_k)}/\log{k}
= (\beta- 1)/\bigl[1-(D-1) (\beta-1)\bigr]. %
\end{equation}
\end{thmm}
Theorem~\ref{thm131} implies that,
when $\bar{F}(s) \sim cs^{-\beta}$ as $s\rightarrow\infty$,
for $\beta\in(1,D/(D-1))$ and $c>0$,
then $P_k = k^{-(1+o(1))\gamma_D(\beta)}$, where $\gamma_D(\beta)$ is
the right-hand side of~(\ref{eq132}). Note that
$\gamma_D(\beta)\searrow0$ as $\beta\searrow1$ and
$\gamma_D(\beta)\nearrow\infty$ as $\beta\nearrow D/(D-1)$.

%
\begin{thmm}
\label{thm141}
Consider a family of JSQ networks as in Theorem~\ref{thm121}, with~(\ref{eq122}) replaced by
%
%
\begin{equation}
\label{eq142} %
c_1 \le\varliminf_{s\rightarrow\infty}s^{D/(D-1)}
\bar{F}(s) \le\varlimsup_{s\rightarrow\infty}s^{D/(D-1)}\bar
{F}(s) \le
c_2 %
\end{equation}
for some $0 < c_1 \le c_2 < \infty$.
Then, for appropriate $r_D(c_2) > 0$ and $s_D(c_1) < \infty$,
%
%
\begin{equation}
\label{eq143} %
r_D(c_2) \le
\varliminf_{k\rightarrow\infty} (1/k)\log{(1/P_k)} \le
\varlimsup_{k\rightarrow\infty} (1/k)\log{(1/P_k)} \le s_D(c_1),
\end{equation}
where
%
%
\begin{eqnarray}
\label{eq144} %
 r_D(c_2) &\nearrow&\infty\qquad
\mbox{as } c_2 \searrow0,
\nonumber
\\[-8pt]
\\[-8pt]
\nonumber
s_D(c_1) &\searrow& 0 \qquad\mbox{as } c_1
\nearrow\infty. %
\end{eqnarray}
\end{thmm}

Theorem~\ref{thm141} implies that when $\bar{F}(s)\sim cs^{-D/(D-1)}$
as $s\rightarrow\infty$, then
$P_k$ decreases exponentially fast in the
sense of (\ref{eq143}). Because of (\ref{eq144}), the exponent
depends strongly on the choice
of $c$.

When $\bar{F}(\cdot)$ satisfies (\ref{eq122}) for a given $\beta>
1$, the asymptotic behavior
of $P_k$ behaves according to (\ref{eq123}) or (\ref{eq132}),
depending on whether
$D > \beta/(\beta- 1)$ or $D < \beta/(\beta- 1)$. In applications
where there is a substantial
penalty for a moderately large number of jobs at a queue (resulting,
e.g., in memory overflow),
it is therefore important to choose $D > \beta/(\beta- 1)$. This
distinction does not occur when
$\bar{F}(\cdot)$ has an exponential tail, since any choice of $D\ge
2$ produces a doubly exponential
tail for $P_k$, as in (\ref{eq111}). (See~\cite{BLP1} for more detail.)

We point out that the proofs of Theorems~\ref{thm121}--\ref{thm141}
only depend on (\ref{eq112}) for
the existence of an equilibrium environment.
Regardless of how the existence of an equilibrium environment is verified,
(\ref{eq112}) will be needed in order to relate the tail behavior of
$P_k$ for the equilibrium environment
to the tail behavior
for the equilibria of the corresponding family of networks as
\mbox{$N\rightarrow\infty$}.

We also note that, although the phrase ``join the shortest queue
network'' is widely
used in the literature, such
systems are not true networks in the sense that, upon the departure of
a job from a queue,
the job leaves the system instead of being able to return to a
different queue.
However, such systems have been extended to the setting of Jackson networks
(see, e.g.,~\cite{MS} and~\cite{SV}).

This article is organized as follows. In Section~\ref{sec2}, we provide basic
background on the properties of the state
space and Markov process that underlie the JSQ networks. We then define
equilibrium environments and formally state the ansatz.
In Sections~\ref{sec3}--\ref{sec5}, we demonstrate Theorems~\ref{thm121},~\ref{thm131}
and~\ref{thm141}, respectively.
Our approach will be to demonstrate lower bounds and then upper bounds
that yield the theorem.
In each case, the lower bounds will be considerably easier to show.


\subsection*{Notation}
For the reader's convenience, we mention here some of the notation in
the paper. We will employ
$C_1,C_2,\ldots$ to denote positive constants whose precise value is
not of
importance to us. For $z\in\mathbb{R}$, $\lfloor z \rfloor$ and
$\lceil z \rceil$ will denote, respectively, the
integer part of $z$ and the smallest integer at least as
large as $z$.

\section{Markov process background, equilibrium environments and the ansatz}\label{sec2}

In this section, we provide a more detailed description of the
construction of
the Markov processes $X^{(N)}(\cdot)$ that underlie the JSQ networks.
We next define the
corresponding cavity process and its equilibrium environment. We then
employ these concepts
to state the ansatz for JSQ networks. Most of this material is included
in Sections~\ref{sec2} and~\ref{sec3} of
Bramson et al.~\cite{BLP2}. (Related material is also given in \cite
{B1} and~\cite{B2}.)

We define the state space $S^{(N)}$ to be the set
%
\begin{equation}
\label{eq211} \bigl(\mathbb{Z} \times\mathbb{R}^2
\bigr)^N.
\end{equation}
The first coordinate $z^n$, $n=1, \ldots, N$, corresponds to the
number of jobs at the $n$th queue;
the second coordinate $u^n$, $u^n \ge0$, is the amount of time the
oldest job there has already been served;
and the last coordinate $s^n$, $s^n > 0$, is the residual service time.
When $z^n = 0$, set the other two coordinates
equal to $0$. The coordinate $u^n$ will not play a role in the
evolution of $X^{(N)}(\cdot)$ here; we retain it
for comparison with~\cite{BLP2}, where it was used to demonstrate
(\ref{eq112}) under decreasing hazard rates.
(We will employ slightly different notation here than in~\cite{BLP2}.)

For given $N^{\prime} \le N$, $S^{(N^{\prime})}$ is the
\emph{projection} of $S^{(N)}$ obtained by restricting $S^{(N)}$ to
the first $N^{\prime}$ queues;
for $x\in S^{(N)}$, $x^{\prime}\in S^{(N^{\prime})}$ is thus obtained
by omitting the coordinates with $n > N^{\prime}$.
One can also define projections of $S^{(N)}$ onto spaces $S^{(N^{\prime})}$
corresponding to other subsets of $\{1,\ldots,N\}$ analogously,
although these are not needed here.

We define the metric $d^{(N)}(\cdot,\cdot)$ on $S^{(N)}$, with
$d^{(N)}(\cdot,\cdot)$
given in terms of $d^{(N),n}(\cdot,\cdot)$ by
$d^{(N)}(\cdot,\cdot) = (1/N)\sum_{n=1}^N d^{(N),n}(\cdot,\cdot)$.
For given $x_1,x_2\in S^{(N)}$, with the coordinates labelled
correspondingly, set
%
%
\begin{equation}
\label{eq212} d^{(N),n}(x_1,x_2)=
\bigl|z_1^{n} - z_2^{n}\bigr| +
\bigl|u_1^{n} - u_2^{n}\bigr| +
\bigl|s_1^{n} - s_2^{n}\bigr|.
\end{equation}
One can check that the metric $d^{(N)}(\cdot,\cdot)$
is separable and locally compact;
more detail is given on page 82 of~\cite{B1}.
We equip $S^{(N)}$ with the standard Borel $\sigma$-algebra inherited
from $d^{(N)}(\cdot,\cdot)$,
which we denote by $\mathscr{S}^{(N)}$.

The Markov process $X^{(N)}(t)$, $t \ge0$, underlying a given model is defined
to be the right continuous process with left limits, taking values $x$
in $S^{(N)}$, whose
evolution is determined by the model together with the assigned service
discipline.
We denote the random values of the coordinates $z^n$, $u^n$ and $s^n$
taken by $X^{(N)}(t)$, by $Z^n(t)$, $U^n(t)$ and $S^n(t)$. Jobs
are allocated service according to the FIFO discipline; during the
period a job is being served,
$U^n(t)$ increases at rate $1$ and $S^n(t)$ decreases at rate $1$.

Along the lines of page 85 of~\cite{B1}, a filtration $(\mathcal
{F}_{t}^{(N)})$, $t\in
[0,\infty]$, can be assigned to $X^{(N)}(\cdot)$ so that
$X^{(N)}(\cdot)$ is
a piecewise-deterministic Markov process, and hence is Borel right.
This implies that $X^{(N)}(\cdot)$
is strong Markov. (We do not otherwise use Borel right.) The reader is
referred to Davis~\cite{D} for more detail.

%
%
%
%

\subsection*{Equilibrium environments and the ansatz}
In order to state the ansatz, we require some terminology.
We denote by
$\mathcal{E}^{(N,N^{\prime})}$ the projection of the equilibrium measure
$\mathcal{E}^{(N)}$ of the $N$-queue system onto the first
$N^{\prime}$ queues. [Since $X^{(N)}(t)$ is exchangeable when
$X^{(N)}(0)$ is, the choice of
queues will not matter.]

We wish to describe the evolution of individual queues for the limiting
process, as $N\rightarrow\infty$. For this, we
construct a strong Markov process $X^{\mathcal{H}}(t)$, $t\ge0$, on $S^{(1)}$.
We will define $X^{\mathcal{H}}(t)$ similarly to $X^{(1)}(t)$, except
that only a fraction of incoming
potential arrivals at the queue is permitted to arrive at the queue,
with the fraction depending on
the current number of jobs there, and with the fraction decreasing as
the number of jobs increases.

We proceed as follows. Let
$\mathcal{H}$ denote a probability measure on $S^{(1)}$, which we
refer to as the \emph{environment} of the
process $X^{\mathcal{H}}(\cdot)$; we refer to $X^{\mathcal{H}}(\cdot
)$ as the associated \emph{cavity process}.
We define $X^{\mathcal{H}}(\cdot)$ so that \emph{potential arrivals}
arrive according to a rate-$D\alpha$ Poisson process. When such a
potential arrival to the queue
occurs at time $t$,
$X^{\mathcal{H}}(t-)$ is compared with the states of $D-1$ independent
random variables,
each with law $\mathcal{H}$; we refer to
these $D-1$ states at a potential arrival as the \emph{comparison states}.
Choosing from among these $D$ states, the job is assigned
to the state with the fewest number of jobs. (In case of a tie, each of
these states is chosen with equal probability.)
If the job has chosen the state $X^{\mathcal{H}}(t-)$ at the queue, it then
immediately joins the queue; otherwise, the job immediately leaves the
system. In either case, the
independent $D-1$ states employed for this purpose are immediately discarded.

We give the following illustrations,
denoting by $Q_k$ the probability that the environment $\mathcal{H}$
has at least
$k$ jobs. For $D=2$, if a potential arrival occurs at time $t$ and
$X^{\mathcal{H}}(t-) = k$, then the
probability that $X^{\mathcal{H}}(t) = k+1$ is $(Q_k + Q_{k+1})/2$,
and so the rate $\alpha_k$ of an
arrival at the queue is $\alpha(Q_k + Q_{k+1})$. For general $D$, in order
for a potential arrival to arrive at the
queue, it is necessary for all of the $D-1$ comparison states used at
that time to be at least $k$, in
which case the probability of selecting the queue is the reciprocal of
the number of states equal to $k$.
This gives the bounds
%
\begin{equation}
\label{eq234new} %
\alpha Q_k^{D-1} \le
\alpha_k \le\alpha DQ_k^{D-1}. %
\end{equation}

We assume that jobs in the cavity process $X^{\mathcal{H}}(\cdot)$
have the same service distribution $F(\cdot)$
as in the queueing network and are served according to the FIFO service
discipline.
The number of jobs in $X^{\mathcal{H}}(t)$ will be denoted by
$Z^{\mathcal{H}}(t)$,
the amount of time the oldest job has already been served by
$U^{\mathcal{H}}(\cdot)$ and the
residual service time by $S^{\mathcal{H}}(t)$; we will employ $x$,
$z$, $u$ and $s$ for the corresponding
terms in the state space.

When a cavity process $X^{\mathcal{H}}(\cdot)$, with environment
$\mathcal{H}$, is stationary with
the equilibrium measure $\mathcal{H}$ [i.e., $X^{\mathcal{H}}(t)$ has
the distribution $\mathcal{H}$ for all $t$],
we say that $\mathcal{H}$ is an \emph{equilibrium environment}. One can
think of an equilibrium environment as being the restriction of an
equilibrium measure for the JSQ network,
viewed at a single queue, when ``the total number of queues $N$ is
infinite.'' More background on
the cavity process is given in~\cite{BLP1}.

We now state the ansatz. Here, $\stackrel{v}{\rightarrow}$ on
$S^{(N^{\prime})}$
denotes convergence in total variation with respect to the metric
$d^{N^{\prime}}(\cdot,\cdot)$ on $S^{(N^{\prime})}$.
\begin{ansatz*}
\label{Ansatz}
Consider a family of JSQ networks, with given $D\ge2$ and
$N=D,D+1,\ldots,$ where the $N$th network
has Poisson rate-$\alpha N$ input, with \mbox{$\alpha< 1$}, and where service
at each queue is FIFO, with
distribution $F(\cdot)$ having mean~$1$. Then, \textup{(a)} for each $N^{\prime}$,
%
%
\begin{equation}
\label{eq231} \mathcal{E}^{(N,N^{\prime})} \stackrel{v}
{\rightarrow} 
\mathcal{E}^{(\infty,N^{\prime})}\qquad \mbox{as } N \rightarrow\infty,
\end{equation}
where $\mathcal{E}^{(\infty,N^{\prime})}$ is the $N^{\prime}$-fold
product of $\mathcal{E}^{(\infty,1)}$.
Moreover, \textup{(b)} $\mathcal{E}^{(\infty,1)}$ is the unique equilibrium
environment associated with this family of networks.
\end{ansatz*}

As was mentioned in the \hyperref[sec1]{Introduction}, this ansatz was
demonstrated in Bramson et al.~\cite{BLP2} when the
service time distribution $F(\cdot)$ has a decreasing hazard rate
$h(\cdot)$
[i.e., $h(s)=F^{\prime}(s)/\bar{F}(s)$ is nonincreasing in $s$] and
for general service distributions when
the arrival rates are small enough.

In order to demonstrate Theorems~\ref{thm121}--\ref{thm141}, we will
analyze the cavity process
$X^{\mathcal{H}}(\cdot)$ with its unique equilibrium environment
$\mathcal{H} = \mathcal{E}^{(\infty,1)}$.
In particular, we will analyze $\mathcal{E}^{(\infty,1)}$ over a
\emph{cycle} starting and ending
at the state $0$. (The state where the number of jobs $z$ is $0$.)
Letting $\nu$ denote the time at
which $X^{\mathcal{H}}(\cdot)$ first returns to~$0$ after visiting
another state, the first cycle is the
random time interval $[0,\nu]$. For any $k\ge1$, we will denote by
$V_k$ the \emph{occupation time}
at states $x$, with $z\ge k$, over $[0,\nu]$, that is,
\[
V_k = \int_{0}^{\nu}\mathbf{1}
\bigl\{Z^{\mathcal{H}}(t)\ge k\bigr\} \,dt.
\]
Setting $m_0 = E[\nu]$, the mean return time to $0$, one has
%
%
\begin{equation}
\label{eq233} %
P_k = m_0^{-1}E[V_k],
\end{equation}
where $P_k$ is the probability there are at least $k$ jobs in the
equilibrium environment.

Letting $\alpha_k$ denote the arrival rate of jobs for $X^{\mathcal
{H}}(\cdot)$ when $z=k$, one has
%
%
\begin{equation}
\label{eq234} %
\alpha P_k^{D-1} \le
\alpha_k \le\alpha DP_k^{D-1}, %
\end{equation}
which is the analog of (\ref{eq234new}).
Since the departure of jobs from the queue is deterministic, being a
function of the residual service time
$s$, (\ref{eq234}) gives a reasonably explicit description of the
transition rates for
$X^{\mathcal{H}}(\cdot)$. Together with (\ref{eq233}), (\ref
{eq234}) will provide the basis
for our demonstration of Theorems~\ref{thm121}--\ref{thm141} and will
be used throughout the paper.

\section{\texorpdfstring{The case where $\beta> D/(D-1)$}
{The case where beta > D/(D-1)}}\label{sec3}

In this section, we demonstrate Theorem~\ref{thm121}; we do this by
demonstrating lower
and upper bounds that are needed for the theorem in Propositions \ref
{prop311} and~\ref{prop361}. Each
of these bounds is expressed in terms of a recursion relation for
$P_k$. In order to obtain Theorem~\ref{thm121}
from these recursions, we employ Proposition~\ref{prop3p11}, which
analyzes such recursions by utilizing a
standard framework involving rational generating functions. The section
is organized as follows. After stating
Propositions~\ref{prop311} and~\ref{prop361}, we state and prove
Proposition~\ref{prop3p11}. We next
employ the three propositions to demonstrate Theorem~\ref{thm121}. We
then provide the relatively quick proof
of Proposition~\ref{prop311} and the longer proof of Proposition~\ref
{prop361}, in the following subsections.

In both propositions, we set $k_1 = \lceil k-\beta\rceil$ (or,
equivalently, $\lfloor\beta\rfloor= k-k_1$)
and $\hat{\beta} = \beta- \lfloor\beta\rfloor$.

%
\begin{prop}
\label{prop311}
Consider a family of JSQ networks, with given $D\ge2$ and
$N=D,D+1,\ldots,$ where the $N$th network
has Poisson rate-$\alpha N$ input, with $\alpha< 1$, and where service\vadjust{\goodbreak}
at each queue is FIFO, with
distribution $F(\cdot)$ having mean $1$. Assume that (\ref{eq112})
holds. Then, for appropriate
$\eee>0$ and all $k$,
%
%
\begin{equation}
\label{eq312} %
P_k \ge(\eee/8k)^k \prod
_{i=0}^{k-1} P_{i}^{D-1}.
\end{equation}
If moreover, for some $s_0 \ge1$,
%
%
\begin{equation}
\label{eq315} %
\bar{F}(s) \ge s^{-\beta} \qquad\mbox{for } s\ge
s_0, %
\end{equation}
with $\beta\in(D/(D-1),\infty)$, then, for appropriate $\eee>0$ and
all $k$,
%
%
\begin{equation}
\label{eq316} %
P_k \ge\eee3^{-k} \Biggl(\prod
_{i=k_1 +1}^{k-1} P_{i}^{D-1}
\Biggr) P_{k_1}^{\hat{\beta} (D-1)}. %
\end{equation}
%
%
\end{prop}

%
%
%
%
%
%
%
%
%
%
%
%
%

%
\begin{prop}
\label{prop361}
Consider a family of JSQ networks, with given $D\ge2$ and
$N=D,D+1,\ldots,$ where the $N$th network
has Poisson rate-$\alpha N$ input, with $\alpha< 1$, and where service
at each queue is FIFO, with
distribution $F(\cdot)$ having mean $1$. Assume that (\ref{eq112})
holds and that, for some $s_0 \ge1$,
%
%
\begin{equation}
\label{eq362} %
\bar{F}(s) \le s^{-\beta} \qquad\mbox{for } s\ge
s_0, %
\end{equation}
with $\beta\in(D/(D-1),\infty)$. If $\beta$ is not an integer,
then, for appropriate $\fff$ and all~$k$,
%
%
\begin{equation}
\label{eq363} %
P_k \le\fff k^{\beta+1} \Biggl(\prod
_{i=k_1 +1}^{k-1} P_{i}^{D-1}
\Biggr) P_{k_1}^{\hat{\beta} (D-1)}. 
\end{equation}
%
If $\beta$ is an integer, then, for each $\delta> 0$, appropriate
$\fff$ and all $k$,
%
%
\begin{equation}
\label{eq364} %
P_k \le\fff k^{\beta+1} \Biggl(\prod
_{i=k_1 +2}^{k-1} P_{i}^{D-1}
\Biggr) P_{k_1+1}^{(1-\delta)(D-1)}. 
\end{equation}
%
%
%
%
%
%
%
\end{prop}

To employ the recursions in (\ref{eq316}) and (\ref{eq363})--(\ref
{eq364}) of
Propositions~\ref{prop311} and~\ref{prop361} in the proof of Theorem
\ref{thm121},
we will analyze the asymptotic behavior of the recursions in (\ref{eq3p12}).
%
%
\begin{prop}
\label{prop3p11}
Suppose that $R_k$ satisfies
%
%
\begin{equation}
\label{eq3p12} %
R_k = (D-1) \Biggl(\sum
_{i = k-\ell+1}^{k-1} R_i + \eta R_{k-\ell}
\Biggr) \qquad\mbox{for } k\ge1, %
\end{equation}
with $R_k = 1$ for $k = -\ell+1,\ldots,-1,0$, where $\ell,D\ge2$
and $\eta\in[0,1]$. Then,
setting $\beta= \ell+ \eta-1$,
%
%
\begin{equation}
\label{eq3p13} %
\lim_{k\rightarrow\infty}\frac{1}{k}
\log_D R_k = q_D (\beta) %
\end{equation}
for some $q_D (\beta) \in(0,1)$. Moreover, $q_D (\beta)$ is
continuous in $\beta$ and
$q_D (\beta) \nearrow1$ exponentially fast as $\beta\nearrow\infty$.\vadjust{\goodbreak}
\end{prop}
\begin{pf}
The recurrence (\ref{eq3p12}) is a special case of linear recursions
of the form
%
%
\begin{equation}
\label{eq3p20} R_k + \sum_{i=1}^{\ell}a_i
R_{k-i} = 0,
\end{equation}
with $a_i \in\mathbb{C}$ and general $R_{-\ell+1},\ldots,R_0$.
It is well known that (see, e.g., Stanley~\cite{ST}, page 202)
%
%
\begin{equation}
\label{eq3p21} %
R_k = \sum_{i=1}^{j}
P_i(k) \gamma_i^k %
\end{equation}
for each $k$, where $\gamma_i$ are distinct, $P_i(k)$ is a
polynomial in $k$ of degree strictly less than $\ell_i$, and
%
%
\begin{equation}
\label{eq3p22} 1 + \sum_{i=1}^{\ell}
a_i x^i = \prod_{i=1}^{j}(1
- \gamma_i x)^{\ell_i},
\end{equation}
with $\sum_{i=1}^{j}\ell_i = \ell$. Moreover the converse holds,
that is,
if (\ref{eq3p21}) and (\ref{eq3p22}) both hold, then so does (\ref{eq3p20}).

For $R_k$ given by (\ref{eq3p12}), it is not difficult to check that
there is
exactly one value~$\gamma_i$, say $\gamma_1$, that is real and
positive, that
$\gamma_1$ varies continuously in $\eta$, and
moreover that $\gamma_1$ satisfies $\gamma_1 > 1$, since $a_i < 0$ and
$\sum_{i=1}^{\ell}a_i < -1$. (Descartes' rule of signs in fact
implies that
$1/\gamma_1$ is a simple root.) Also, because $a_i < 0$, and possesses
both odd and even
indices, $|\gamma_i| < \gamma_1$ for $i\neq1$. Since the initial
data given below
(\ref{eq3p12}) are all positive, any solution of (\ref{eq3p12}) is
majorized by
this particular solution, up to a multiplicative constant; so,
$P_1(\cdot)\not\equiv0$.
The limit in (\ref{eq3p13}), with $q_D(\beta) = \log_D\gamma_1 >
0$, follows from
these observations.

We still need to examine the limiting behavior of $q_D(\beta)$ as
$\beta\rightarrow\infty$.
Dividing both sides in (\ref{eq3p12}) by $R_k$, then substituting
(\ref{eq3p21})
for each of the terms, and letting $k\rightarrow\infty$ implies that
\begin{eqnarray*}
1 &=& (D-1) \bigl(x + x^2 + \cdots+ x^{\ell-1} +
\eta x^{\ell} \bigr)
\\
&=& (D-1) \bigl(x - (1-\eta)x^{\ell} - \eta x^{\ell+ 1} \bigr)/(1-x)
\end{eqnarray*}
for $x = 1/\gamma_1 = D^{-q_D(\beta)}$. This again uses $\gamma_1 >
|\gamma_i|$
for $i\neq1$. Hence,
%
%
\begin{equation}
\label{eq3p23} %
Dx-1 = (D-1) \bigl((1-\eta)x^{\ell} + \eta
x^{\ell+ 1} \bigr). %
\end{equation}
Note that $x\in(0,1)$ and that, since $q_D(\beta)$ is increasing in~$\beta$,
$x$ is decreasing in~$\beta$. Since the right-hand side goes to $0$
exponentially fast as $\ell\nearrow\infty$, and hence as $\beta
\nearrow\infty$,
it follows that $x\searrow1/D$ exponentially fast as $\beta\nearrow
\infty$, which also implies
$q_D(\beta)\nearrow1$ exponentially fast, as desired. Note that the precise
exponential rate of convergence can be obtained by inserting this limit back
into the right-hand side of (\ref{eq3p23}).\vadjust{\goodbreak}
\end{pf}

Applying Proposition~\ref{prop3p11} to Propositions~\ref{prop311} and
\ref{prop361}, we
now demonstrate Theorem~\ref{thm121}.
\begin{pf*}{Proof of Theorem~\ref{thm121}}
Setting $Q_k = e^{R_k}$, where $R_k$ is given in (\ref{eq3p12}), one has
%
%
\begin{equation}
\label{eq3q11} %
Q_k = \Biggl(\prod
_{i=k-\ell+1}^{k-1} Q_{i}^{D-1} \Biggr)
Q_{k-\ell}^{\eta(D-1)}, %
\end{equation}
with $Q_k = e$ for $k=-\ell+1,\ldots,-1,0$. We proceed to compare
$Q_k$ with
$1/P_k$, where $P_k$ satisfies one of (\ref{eq316}), (\ref{eq363})
and (\ref{eq364}).

Comparison of $Q_k$ with $1/P_k$, with $\eta= \hat{\beta}$, $\ell=
\lfloor\beta\rfloor= k - k_1$ and
$P_k$ satisfying~(\ref{eq316}), provides an upper bound on the limit
in (\ref{eq123}).
To see this, we first set $\tilde{Q}_k = M^{-k}Q_k$, for given $M>1$.
Since $(D-1)(\beta-1)>1$, by
substituting into~(\ref{eq3q11}),
one can check that, for large enough $M$ and $k$,
%
%
\begin{equation}
\label{eq3q12} %
\tilde{Q}_k \ge\hhh b^k
\Biggl(\prod_{i=k-\ell+1}^{k-1}
\tilde{Q}_{i}^{D-1} \Biggr) \tilde{Q}_{k-\ell}^{\eta(D-1)}
\end{equation}
for any fixed choice of $\hhh$ and $b$, in particular, for $\hhh=
1/\eee$ and $b=3$, where $\eee$
is chosen as in Proposition~\ref{prop311}. Moreover, on account of
(\ref{eq3p13}),
%
%
\begin{equation}
\label{eq3q13} %
\lim_{k\rightarrow\infty}(1/k)\log_D{\log{(
\tilde{Q}_k)}} = q_D(\beta), %
\end{equation}
where, in particular, $q_D(\beta) > 0$, and hence $\tilde{Q}_k
\rightarrow\infty$ as $k\rightarrow\infty$.

We observe that $1/P_k$ satisfies the inequality that is analogous to
that for $P_k$ in (\ref{eq316}), but
with the inequality reversed and prefactors $3^k/\eee$ instead of
$\eee/3^k$.
Comparing $\tilde{Q}_k$ with $1/P_k$ therefore implies that, for large
enough $n$ not depending on $k$,
\[
1/P_k \le\tilde{Q}_{k+n}.
\]
The upper bound for (\ref{eq123}) therefore follows from (\ref
{eq3q13}) for the same choice of
$q_D(\beta)$, which we recall is continuous in $\beta$.
The limit in (\ref{eq124}) also follows from Proposition~\ref{prop3p11}.

Comparison of $Q_k$ with $1/P_k$ also provides a lower bound on the
limit in (\ref{eq123}). In the
case where $\beta$ is nonintegral, we choose $\eta$ and $\ell$ as
before, with
$\eta= \hat{\beta}$, $\ell= \lfloor\beta\rfloor= k - k_1$; note that
$P_k$ satisfies the upper bound in (\ref{eq363}). We
proceed as in the first part, but instead set $\tilde{Q}_k =
M^{k}Q_k$, for given $M>1$.
One can check that, for large enough $M$ and $k$,
%
%
\begin{equation}
\label{eq3q14} %
\tilde{Q}_k \le\hhh b^k
\Biggl(\prod_{i=k-\ell+1}^{k-1} Q_{i}^{D-1}
\Biggr) Q_{k-\ell}^{\eta(D-1)} %
\end{equation}
for any choice of $\hhh>0$ and $b>0$. As before, (\ref{eq3q13}) holds.

The terms $1/P_k$ satisfy the inequality that is the analog of (\ref
{eq363}). Also,
$1/P_k \rightarrow\infty$ as $k\rightarrow\infty$.
Comparing\vadjust{\goodbreak} $\tilde{Q}_k$ with $1/P_k$ therefore implies that, for large
enough $n$ not depending on $k$,
%
%
\begin{equation}
\label{eq3q15} 1/P_{k+n} \ge\tilde{Q}_{k}.
\end{equation}
The lower bound for (\ref{eq123}) therefore follows from (\ref
{eq3q13}) when $\beta$ is nonintegral.

The reasoning in the case where $\beta$ is integral is similar, but
with the difference that we now choose
$\eta= 1 - \delta$, $\ell= \beta-1 = k - k_1 -1$, where $\delta\in
(0,1)$ is arbitrary. Now, $P_k$ satisfies
the upper bound in (\ref{eq364}). We proceed as in the nonintegral
case, once again obtaining
(\ref{eq3q14}). Comparing $1/P_k$ with $\tilde{Q}_k$ again produces
(\ref{eq3q13}), except that the limit
is now $q_D(\beta- \delta)$ because of our choice of $\eta$. By
Proposition~\ref{prop3p11},
$q_D(\cdot)$ is continuous in its argument. Therefore, letting $\delta
\searrow0$ produces the same limit
as in the nonintegral case, and hence implies the lower bound for (\ref
{eq123}) in the case where $\beta$ is integral.

We still need to demonstrate that when (\ref{eq122}) holds with $\beta
= \infty$, then (\ref{eq123})
holds with $q_D(\infty) =1$. The lower bound in (\ref{eq123}) holds
on account of (\ref{eq124}). The upper
bound is not difficult to show and does not require Proposition~\ref
{prop3p11}; we proceed to show the bound.

We will show by induction that, for all $k$,
%
%
\begin{equation}
\label{eq3q16} P_k \ge(\eee/8k)^{k D^k}, 
\end{equation}
where $\eee$ is as chosen as in (\ref{eq312}), which we assume WLOG
is at most $1$.
To see~(\ref{eq3q16}), note that if it holds for all $i=0,\ldots,k-1$
then this,
together with~(\ref{eq312}), implies that
\begin{eqnarray*}
P_k &\ge&(\eee/8k)^k \prod
_{i=0}^{k-1} \bigl[(\eee/8i)^{iD^{i}}
\bigr]^{D-1}
\\
&\ge&(\eee/8k)^{(k-1)(D^k -1) +k} \ge(\eee/8k)^{k D^k}. %
\end{eqnarray*}
The upper bound in (\ref{eq123}), with $q_D(\infty) = 1$, follows
immediately from (\ref{eq3q16}).
\end{pf*}

\subsection*{Demonstration of Proposition \protect\ref{prop311}}

The proof of Proposition~\ref{prop311} is quick. To obtain the lower
bounds in both (\ref{eq312})
and (\ref{eq316}), it suffices to construct a path along which
$Z^{\mathcal{H}}(t)$ increases
from $0$ to $k$ within the first cycle. This is done, in both cases, by
allocating the same amount of
time to each of the first $k$ arrivals, which are also required to
occur before the first departure.

\begin{pf*}{Proof of Proposition~\ref{prop311}}
Consider the cavity process $X^{\mathcal{H}}(\cdot)$ with
$X^{\mathcal{H}}(0) = 0$.
In order to show (\ref{eq312}) and (\ref{eq316}), we obtain lower
bounds on the
expected amount of time $E[V_k]$ over which $Z^{\mathcal{H}}(t) \ge k$
before $X^{\mathcal{H}}(\cdot)$ returns to $0$. We first show (\ref{eq312}).

We consider the event $A$ where the first service time $S$ is at least
$1/2$ and
the first $k$ arrivals occur by time $1/4$. The latter event contains
the event
where each of the first $k$ arrivals occurs not more than $1/4k$ units
of time after the
previous arrival, starting at time $0$.

Conditioned on there being $i$ jobs in the queue, jobs arrive at rate
$\alpha_i \ge\alpha P_{i}^{D-1}$,
and so the probability of such an arrival occurring over an interval of
length $1/4k$ is at least
$1 - \exp\{-\alpha P_{i}^{D-1}/4k\}$. So, given that $S\ge1/2$, the
probability that all $k$ of
these arrivals occur by time $1/4$ is at least
%
%
\begin{equation}
\label{eq371} %
\prod_{i=0}^{k-1}
\bigl(1 - \exp\bigl\{-\alpha P_{i}^{D-1}/4k\bigr\} \bigr).
\end{equation}

The event $S \ge1/2$ occurs with some positive probability $c$
depending on $F(\cdot)$ and,
under the event $A$, the departure time
for the first job occurs at least $1/4$ after the last of the first $k$
arrivals.
So, the expected amount of time in $[1/4,1/2]$, during which
$Z^{\mathcal{H}}(t) \ge k$ and
before $X^{\mathcal{H}}(\cdot)$ has returned to $0$, is at least
%
%
\begin{equation}
\label{eq372} %
\frac{c}{4}\prod_{i=0}^{k-1}
\bigl(1 - \exp\bigl\{-\alpha P_{i}^{D-1}/4k\bigr\} \bigr),
\end{equation}
which is therefore a lower bound for $E[V_k]$. It therefore follows
from (\ref{eq233}) that
%
%
\begin{equation}
\label{eq373} 
P_k \ge\frac{c}{4m_0}\prod
_{i=0}^{k-1} \bigl(1 - \exp\bigl\{-\alpha
P_{i}^{D-1}/4k\bigr\} \bigr) \ge\frac{c}{4}(\alpha/
8k)^k \prod_{i=0}^{k-1}
P_{i}^{D-1}, %
\end{equation}
which implies (\ref{eq312}) for appropriate $\eee$.

We next show (\ref{eq316}) under the assumption (\ref{eq315}). For
this, we set
%
\begin{equation}
\label{eq374} s_1 = 2k/\bigl(\alpha P_{k_1}^{D-1}
\bigr).
\end{equation}
One can reason analogously as through (\ref{eq372}), but by replacing
the time interval $[0,1/2]$ by
$[0,s_1]$ and employing $s_1/2k$ for the allotted time for each of the
$k$ arrivals. One obtains that
the expected amount of time in $[s_1/2,s_1]$, during which $Z^{\mathcal
{H}}(t) \ge k$ and
before $X^{\mathcal{H}}(\cdot)$ has returned to $0$, is at least
%
%
\begin{equation}
\label{eq381} %
\frac{s_1}{2} \bar{F}(s_1) \prod
_{i=0}^{k-1} \bigl(1 - \exp\bigl\{-\alpha
s_1 P_{i}^{D-1}/2k\bigr\} \bigr). %
\end{equation}

Choose $k$ large enough so that $s_1 \ge s_0$, where $s_0$ is as in~(\ref{eq315}) and $s_1$
is as in~(\ref{eq374}). Since $e^{-x} \le(1-x/2)\vee1/2$ for $x\ge
0$, this is at least
\begin{eqnarray*}
&&2^{-(k_1+2)}s_1^{-(\beta- 1)}(\alpha
s_1/4k)^{k-k_1-1} \prod_{i=k_1+1}^{k-1}P_{i}^{D-1}
\\
&&\qquad \ge2^{-k}(\alpha/4k)^{\beta} \Biggl(\prod
_{i=k_1+1}^{k-1}P_{i}^{D-1} \Biggr)
P_{k_1}^{\hat{\beta}(D-1)}, %
\end{eqnarray*}
where the inequality follows from (\ref{eq374}) and $k-k_1 = \beta-
\hat{\beta}$. Consequently,
\[
E[V_k] \ge2^{-k}(\alpha/4k)^{\beta} \Biggl(\prod
_{i=k_1+1}^{k-1}P_{i}^{D-1}
\Biggr) P_{k_1}^{\hat{\beta}(D-1)}.
\]
Again applying (\ref{eq233}), it follows that, for large enough $k$
(depending on $\alpha$ and $\beta$),
\[
P_k \ge3^{-k} \Biggl(\prod_{i=k_1+1}^{k-1}P_{i}^{D-1}
\Biggr) P_{k_1}^{\hat{\beta}(D-1)},
\]
which implies (\ref{eq316}).
\end{pf*}

\subsection*{Demonstration of Proposition \protect\ref{prop361}}

In order to demonstrate Proposi-\break tion~\ref{prop361}, we will employ
Lemma~\ref{lem391} below; the lemma
will also be employed in the demonstration of Propositions \ref
{prop441} and~\ref{prop531}. (A substantially
more intricate variant of the proof of Lemma~\ref{lem391} will be
needed for the proof of Proposition~\ref{prop461}.)
Lemma~\ref{lem391} provides upper bounds involving $R(k,s)$,
$H(n)$ and $\rho(k,s)$, for $k\ge1$, $s\ge0$ and $n\ge0$, which are
defined as follows.

For $s > 0$, $R(k,s)$ is the expected return time of the cavity
process $X^{\mathcal{H}}(\cdot)$ (with equilibrium environment
$\mathcal{H}$)
to the empty state $0$, from $X^{\mathcal{H}}(0)$ with $Z^{\mathcal
{H}}(0)=k$ and
$S^{\mathcal{H}}(0) =s$. We set $R(k,0) = \lim_{s\searrow0}R(k,s)$,
which is also the expected return time to
$0$ just after departure of a job, but without knowledge of the
residual service time of the job that is beginning
service. The quantity $H(n)$ is the
number of jobs, for this process, at the time when the $(n+1)$st job
has just departed, for example,
$H(0)$ is the number of jobs just after departure of the job originally
in service. The stopping time $\rho(k,s)$ is the
first time $n$ at which $H(n) =0$.

We also denote by $Y_n$ the service time of the $(n+1)$st job
(with $Y_0=s$ being the service time of the job originally in service),
and set
$T_{\ell} = \sum_{n=0}^{\ell}Y_n = \sum_{n=1}^{\ell}Y_n + s$. Note
that $Y_1,Y_2,\ldots$ are i.i.d. with distribution
function $F(\cdot)$, which, as always, is assumed to have mean $1$.

%
%
\begin{lem}
\label{lem391}
Let $R(\cdot,\cdot)$ and $\rho(\cdot,\cdot)$ be defined as above.
Then, for large enough $N_0$,
%
%
\begin{equation}
\label{eq392} %
R(k,s) \le2(k+s+N_0) %
\end{equation}
and
%
%
\begin{equation}
\label{eq393} %
E\bigl[\rho(k,s)\bigr] \le2(k+s/2 +
N_0) %
\end{equation}
for all $k$ and $s$.
\end{lem}
\begin{pf}
It is not difficult to see that (\ref{eq392}) follows from (\ref{eq393}).
By applying Wald's equation to $T(\cdot)$ and $\rho(\cdot,\cdot)$\vadjust{\goodbreak}
(with respect to the underlying $\sigma$-algebra generated by
$X^{\mathcal{H}}(\cdot)$), one obtains
%
%
\begin{eqnarray*}
R(k,s) = E[T_{\rho(k,s)}] = E \Biggl[\sum
_{n=1}^{\rho(k,s)}Y_n \Biggr] + s = E\bigl[
\rho(k,s)\bigr]E[Y_1] + s \le2(k+s+N_0), %
\end{eqnarray*}
with the inequality following from (\ref{eq393}) and $E[Y_1] = 1$.

In order to show (\ref{eq393}), we consider the process
%
%
\begin{equation}
\label{eq3101} %
M(n) = H(n) + n/2 - N_1\exp\bigl\{-
\theta\bigl(H(n) \wedge k_0\bigr)\bigr\}. %
\end{equation}
For appropriate choices of $N_1,\theta>0$ and $k_0\in\mathbb{Z}_+$,
we claim $M(n)$ is a supermartingale, with respect to the filtration
$\mathcal{G}_n = \sigma(H(0),\ldots, H(n))$, after restricting to times
$n$, with $n\le\rho(k,s)$, and then stopping the process.

These three constants are chosen as follows. We choose $k_0$ large
enough so that $\alpha D P_{k_0 +1}^{D-1} \le1/2$.
For $H(n)> k_0$, one can check that the supermartingale inequality
%
%
\begin{equation}
\label{eq3102} %
E\bigl[M(n+1)| \mathcal{G}_n\bigr] \le
M(n) %
\end{equation}
is satisfied---the arrival rate of jobs is at most $1/2$ over the time
interval $(T_{n-1}, T_{n}]$ during which the $(n+1)$st job is served,
which has mean length $1$, and so
\[
E\bigl[H(n+1)| \mathcal{G}_n\bigr] \le H(n) - 1/2.
\]
%

In order to analyze $M(n+1)$ when $H(n)\le k_0$, we set
\[
M_1(n) = -\exp\bigl\{-\theta\bigl(H(n) \wedge k_0
\bigr)\bigr\}.
\]
%
We choose $\theta$ large enough so that, for
some $\varepsilon> 0$ and all $H(n) \le k_0$,
%
%
\begin{equation}
\label{eq3103} %
E\bigl[M_1(n+1)| \mathcal{G}_n
\bigr] \le M_1(n) - \varepsilon. %
\end{equation}
This requires a standard computation using the convexity of the
exponential function and the upper
bound $\alpha D$ on the arrival rate of jobs. [Since $H(\cdot)$ may
have positive drift,
$\theta$ may need to be chosen large.]

We also choose $N_1$ so that $\varepsilon N_1 \ge\alpha D + 1/2$.
Together with (\ref{eq3103}),
this implies (\ref{eq3102}) also holds for $H(n)\le k_0$.
Consequently, $M(n)$ is a supermartingale,
as claimed.

In order to demonstrate (\ref{eq393}), we will apply the optional
sampling theorem to $M(\cdot)$
stopped at times $\rho_n(k,s) = \rho(k,s) \wedge n$. First note that
%
%
\begin{equation}
\label{eq3104} %
E\bigl[M(0)\bigr] \le E\bigl[H(0)\bigr] \le k + s/2
\end{equation}
for $k\ge k_0$, since the arrival rate of jobs is bounded above by
$1/2$. Also,
for given $s$, $E[H(0)]$ is increasing
as a function of $k$, the number of jobs in the cavity process at time $0$.
Together with (\ref{eq3104}), this implies that, for all $k$,
%
%
\begin{equation}
\label{eq3105} %
E\bigl[M(0)\bigr] \le(k\vee k_0) + s/2
\le k + s/2 + k_0. %
\end{equation}

Since the supermartingale $M(\cdot)$ is bounded from below,
application of the optional sampling
theorem to $\rho_n(k,s)$ implies that
\[
E\bigl[M\bigl(\rho_n(k,s)\bigr)\bigr] \le E\bigl[M(0)\bigr] \le
k +
s/2 + k_0,\vadjust{\goodbreak}
\]
and hence
\[
0\le E\bigl[H\bigl(\rho_n(k,s)\bigr)\bigr] \le k + s/2 +
k_0 + N_1 - E\bigl[\rho_n(k,s)\bigr]/2.
\]
Solving for $E[\rho_n(k,s)]$ implies
\[
E\bigl[\rho_n(k,s)\bigr] \le2(k + s/2 + k_0 +
N_1) = 2(k + s/2 + N_0)
\]
for $N_0 = k_0 + N_1$. Letting $n\rightarrow\infty$ implies (\ref{eq393}).
\end{pf}

Lemma~\ref{lem391} provides an upper bound on the expected time over a cycle
during which there are at least $k$ jobs, provided such a state has
already been attained.
Below, we will obtain an upper bound on the probability of attaining
such a state and combine this with (\ref{eq392}).

In order for $X^{\mathcal{H}}(\cdot)$, starting at $0$, to attain a
state with $k$ jobs, it must first attain states
with $k_1+1,k_1+2,\ldots,k-1$ jobs, where $k_1$ has been specified in
the previous subsection. (It turns
out that including states with fewer jobs
in this sequence will not improve our bounds.) We let $\sigma
_{k_1+1},\ldots,\sigma_{k}$ denote
the number of jobs that have already departed when such a state is
first attained [e.g.,
$\sigma_i = 0$ means that the first job is still being served at the
time $t$ when $Z^{\mathcal{H}}(t) = i$ first
occurs].

One trivially has
\[
0\le\sigma_{k_1+1}\le\sigma_{k_1+2}\le\ldots\le
\sigma_k.
\]
Partition $\{k_1+1,k_1+2,\ldots,k\}$ so that $i \neq i^{\prime}$ are
in the same subset if
$\sigma_i = \sigma_{i^{\prime}}$, that is, the times $t_i$ and
$t_{i^{\prime}}$ at which
$Z^{\mathcal{H}}(t_i)=i$ and $Z^{\mathcal{H}}(t_{i^{\prime
}})=i^{\prime}$ first occur are in the same
service time interval. One can write such a partition as
%
%
\begin{equation}
\label{eq3111} %
\| i_0+1,\ldots,i_1 \|
i_1+1,\ldots,i_2 \| \ldots\| i_{m-1}+1,
\ldots,i_m \|, %
\end{equation}
with $i_0 = k_1$ and $i_m = k$, when the partition consists of $m$
sets (where $m$ is random). We denote by $\Pi_k$ the set of
all such partitions and by $\pi\in\Pi_k$ an element in the set,
with the notation $i_0(\pi),i_1(\pi),\ldots,i_m(\pi)$
being used when convenient. We will say that a partition $\pi$ occurs
during a cycle when the corresponding
sequence of events occurs, and denote by $A_{\pi}$ the event
associated with the partition.

For each of the sets in (\ref{eq3111}) except the last, there is a
corresponding service interval,
$[T_{n_{\ell}-1}, T_{n_{\ell}})$, with $\ell= 1,\ldots,m-1$, at the
beginning of which
there are strictly less than $i_{\ell-1}$ jobs and at the end exactly
$i_{\ell}$ jobs.
(Since such an interval ends with a departure, the number of jobs at
the beginning of the next service
interval must be one less, which requires the cavity process to
``retrace some of its steps'' before
the number of jobs reaches $i_{\ell}$ again.) For
$\ell=m$, there may be strictly more than $k$ jobs at $T_{n_{\ell}}$;
instead, we consider the restricted interval
$[T_{n_{m}-1}, \tau_k]$, where $\tau_k$ is the first time at which
there are at least $k$ jobs. Unlike
at the end of the other intervals $[T_{n_{\ell}-1}, T_{n_{\ell}})$,
the residual service time $s$ will not be $0$. When $s$ is large, this
will increase the occupation time where $Z^{\mathcal{H}}(t)\ge k$,
which will require us to exercise some
care with our computations.\vadjust{\goodbreak}

Since $k-k_1 \le\beta$, the number of distinct partitions in (\ref
{eq3111}) is at most $2^{\beta}$. In
Proposition~\ref{prop3121} below, we compute an upper bound on $P_k$
using an upper bound on the
expected occupation time corresponding to each partition,
and then by multiplying by $2^{\beta}$. The upper bound in (\ref
{eq3122}) includes
a factor $k^{\beta}$ obtained by employing Lemma~\ref{lem391} repeatedly.
The form of the bounds in (\ref{eq3122}) and (\ref{eq3125}) varies in
different ranges of $s$;
we will therefore find it useful to employ the notation
%
%
\begin{equation}
\label{eq3113} 
%
%
L_{\ell}(s) =
\prod_{i=i_{\ell- 1}}^{i_{\ell} -1 } \bigl[\bigl(\alpha
DP_{i}^{D-1}s\bigr)\wedge1 \bigr]. %
\end{equation}
[$L_{\ell}(\cdot)$ implicitly depends on the partition $\pi$ through
$i_{\ell-1}$ and $i_{\ell}$.]
We will employ $L(s)$ when $i$ goes from $k_1$ to $k-1$, which
corresponds to the trivial partition
in~(\ref{eq3111}) consisting of a single set.

In the proof of Proposition~\ref{prop3121}, we will use the following
elementary Chebyshev integral inequality,
which states that, if $f(s)$ and $g(s)$ are both integrable functions
that are increasing in $s$, then, for any
distribution function $F(\cdot)$,
%
%
\begin{equation}
\label{eq3115} %
\int_{-\infty}^{\infty}f(s)g(s)
F(ds) \ge\int_{-\infty}^{\infty}f(s) F(ds) \cdot\int
_{-\infty}^{\infty}g(s) F(ds). %
\end{equation}
%

\begin{prop}
\label{prop3121}
Consider a family of JSQ networks, with the same assumptions holding as
in Proposition~\ref{prop361},
except that (\ref{eq362}) is not assumed. Then,
for large enough $k$,
%
%
\begin{equation}
\label{eq3122} %
P_k \le3m_{0}^{-1}(6k)^{\beta}
\int_0^{\infty}(k+s)L(s) F(ds). %
\end{equation}
\end{prop}
\begin{pf}
We first claim that the probability of the cavity process
$X^{\mathcal{H}}(\cdot)$, with $Z^{\mathcal{H}}(0)\le i_{\ell-1}$
and $S^{\mathcal{H}}(0) = s$,
attaining $i_{\ell}$ jobs before time $s$ is at most
%
%
\begin{eqnarray}
\label{eq3125} %
\prod_{i=i_{\ell- 1}}^{i_{\ell} - 1}
\bigl(1-\exp\bigl\{-\alpha DP_{i}^{D-1}s\bigr\} \bigr) &\le&
\prod_{i=i_{\ell- 1}}^{i_{\ell} - 1} \bigl[\bigl(\alpha
DP_{i}^{D-1}s\bigr)\wedge1 \bigr]
\nonumber
\\[-8pt]
\\[-8pt]
\nonumber
&=& L_{\ell}(s).
\end{eqnarray}
Under this event, arrivals must occur sequentially over $[0,s]$ at
times $t_i$ when $Z^{\mathcal{H}}(t_i-)=i$, for
$i = i_{\ell-1},\ldots,i_{\ell}-1$, and the rate of such arrivals is
at most $\alpha DP_{i}^{D-1}$.
Since there is at most time $s$ for each arrival, multiplying the
corresponding upper bounds on
the probability of an arrival at each step gives the first bound in
(\ref{eq3125}). The following
inequality is then obtained by applying the inequality $1-e^{-x} \le x
\wedge1$.

Recall that $V_k$ denotes the occupation time over a cycle when
\mbox{$Z^{\mathcal{H}}(t)\ge k$}.
In order for $V_k > 0$, the event $A_{\pi}$ must occur
for some $\pi\in\Pi_k$; hence
$E[V_k] = \sum_{\pi\in\Pi_k} E[V_k; A_{\pi}]$.
We claim that, for any partition $\pi\in\Pi_k$ and large enough~$k$,
%
%
\begin{eqnarray}
\label{eq3126} %
&&E[V_k; A_{\pi}]
\nonumber
\\[-8pt]
\\[-8pt]
\nonumber
&&\qquad \le
(3k)^{m_{\pi}} \prod_{\ell=1}^{m_{\pi}-1}
\biggl(\int_{0}^{\infty}L_{\ell}(s) F(ds)
\biggr) \cdot 3\int_{0}^{\infty}(k+s)L_{m_{\pi}} (s)
F(ds). \hspace*{-35pt}
\end{eqnarray}

To obtain (\ref{eq3126}), we argue by induction, applying (\ref
{eq3125}) at each step.
It suffices to show that, for each step with $\ell< m_{\pi}$, one
obtains an additional
factor $3i_{\ell-1}\int_{0}^{\infty}L_{\ell}(s) F(ds)$ and, for
$\ell= m_{\pi}$, one obtains
the factor $9(i_{m_{\pi} -1})\int_{0}^{\infty}(k+s)L_{m_{\pi}}(s) F(ds)$.
For $\ell\ge2$, the factor $3i_{\ell-1}$ is obtained by applying
(\ref{eq393}), with $s=0$, which
gives an upper bound on the expected number of service intervals
occurring over the remainder of the
cycle, after the service interval corresponding to the $(\ell-1)$st
step ends; also, $i_0 \ge m_0$,
which equals the expected number of service intervals at the beginning
of the cycle.
The other factor is obtained from~(\ref{eq3125}) by integrating
against $F(\cdot)$ and,
for $\ell= m_{\pi}$, by employing (\ref{eq392}) to provide an upper bound
on the expected occupation time $V_k$, again employing (\ref{eq3125})
and then integrating against
$F(\cdot)$.


On the other hand, by repeatedly applying the Chebyshev integral
inequality~(\ref{eq3115})
to (\ref{eq3126}), it follows that, for an arbitrary partition in
(\ref{eq3111}), (\ref{eq3126})
is maximized for the trivial partition. That is, for any partition $\pi
\in\Pi_k$,
the quantity in (\ref{eq3126}) is bounded above by
%
%
\begin{equation}
\label{eq3127} 3(3k)^{\beta} \int_{0}^{\infty}(k+s)L(s)
F(ds).
\end{equation}
%

Since $|\Pi_k| \le2^{\beta}$, it follows from (\ref{eq3126}) and
(\ref{eq3127}) that
\begin{eqnarray*}
P_k &=& m_0^{-1}
E[V_k] = m_0^{-1}\sum
_{\pi\in\Pi_k}E[V_k; A_{\pi}]
\\
&\le& 3m_0^{-1} (6k)^{\beta} \int
_{0}^{\infty}(k+s)L(s) F(ds),
\end{eqnarray*}
%
%
which implies (\ref{eq3122})
\end{pf}

We now complete the proof of Proposition~\ref{prop361}.
\begin{pf*}{Proof of Proposition~\ref{prop361}}
We employ the upper bound for $P_k$ given by~(\ref{eq3122}) for large
enough $k$.
The integral in (\ref{eq3122}) is bounded above by
%
%
\begin{eqnarray}
\label{eq3131} &&2ks_0 \int_{0}^{s_0}L(s)
F(ds) + 2k \int_{s_0}^{\infty}sL(s) F(ds)
\nonumber
\\[-8pt]
\\[-8pt]
\nonumber
&&\qquad \le2\beta\bigl(s_0^{\beta+1}+1\bigr)k \int
_{1}^{\infty}s^{-\beta}L(s) \,ds
\end{eqnarray}
by integrating by parts and absorbing the first term into the second;
note that $L(s)$ is increasing in $s$
on account of (\ref{eq3113}).
We decompose this last integral using intervals of the form
$[1/\alpha DP_{k-1}^{D-1},\infty)$, $[1/\alpha DP_{i-1}^{D-1},
1/\alpha DP_{i}^{D-1})$, for
$i=k_1+1,\ldots,k-1$, and $[1,1/\alpha DP_{k_1}^{D-1})$; we need to
consider the cases where
$\beta$ is and is not an integer separately.

Suppose that $\beta$ is not an integer. Applying (\ref{eq3113}) to
the above integral over
$[1/\alpha DP_{k-1}^{D-1},\infty)$, one has the upper bound
%
%
\begin{equation}
\label{eq3141} %
\int_{1/\alpha DP_{k-1}^{D-1}}^{\infty}
s^{-\beta} \,ds = \frac{1}{\beta-1}(\alpha DP_{k-1})^{(D-1)(\beta-1)}.
\end{equation}
For $i=k_1+1,\ldots,k-1$, one has, over $[1/\alpha DP_{i-1}^{D-1},
1/\alpha DP_{i}^{D-1})$, the upper bounds
%
%
\begin{eqnarray}
\label{eq3142} %
&&\int_{1/\alpha DP_{i-1}^{D-1}}^{1/\alpha DP_{i}^{D-1}} (
\alpha Ds)^{k-i}(P_{k-1}\cdots P_i)^{D-1}s^{-\beta}
\,ds
\nonumber
\\[-8pt]
\\[-8pt]
\nonumber
&&\qquad \le\frac{(\alpha D)^{\beta-1}}{\beta+i-k-1} \bigl(P_{k-1}\cdots P_i
P_{i-1}^{\beta+i-k-1} \bigr)^{D-1}. %
\end{eqnarray}
For the last interval $[1,1/\alpha DP_{k_1}^{D-1})$, one has the upper bound
%
%
\begin{eqnarray}
\label{eq3143} %
&&\int_{1}^{1/\alpha DP_{k_1}^{D-1}} (
\alpha Ds)^{k-k_1}(P_{k-1}\cdots P_{k_1})^{D-1}s^{-\beta}
\,ds
\nonumber
\\[-8pt]
\\[-8pt]
\nonumber
&&\qquad \le\frac{(\alpha D)^{\beta-1}}{1-\hat{\beta}} \bigl
(P_{k-1}\cdots P_{k_1+1}
P_{k_1}^{\hat{\beta}} \bigr)^{D-1}, %
\end{eqnarray}
where we recall that $\hat{\beta} = \beta-k+k_1$. Note that the
lower limits of integration supply the
dominant term in (\ref{eq3141}) and (\ref{eq3142}), whereas the upper
limit supplies the dominant term
in (\ref{eq3143}), because of the choice of $k_1$.

Since $P_i$ is decreasing in $i$, if one ignores the coefficients not
involving powers of $P_i$ on the right-hand sides of (\ref
{eq3141})--(\ref{eq3143}), the largest bounds in (\ref
{eq3141})--(\ref{eq3143})
are given in (\ref{eq3142}), with $i=k_1+1$, and in (\ref{eq3143}),
in each case by the powers of~$P_i$,
%
%
\begin{equation}
\label{eq3144} %
\bigl(P_{k-1}\cdots P_{k_1}^{\hat{\beta}}
\bigr)^{D-1}. %
\end{equation}
The coefficients of these powers are bounded above by terms not
involving~$k$. Employing (\ref{eq3122})
of Proposition~\ref{prop3121}, together with (\ref{eq3131}), one
obtains the bound~(\ref{eq363})
for $P_k$, for appropriate $\fff$ and all $k$.

When $\beta$ is an integer, the computations are similar. The
inequalities in (\ref{eq3141}) and (\ref{eq3143})
are the same as before, as are all of the cases in (\ref{eq3142})
except for $i=k_1+1$. Rather than (\ref{eq3142}),
one obtains the following inequality when $i=k_1+1$:
%
%
\begin{eqnarray}
\label{eq3145} %
&&\int_{1/\alpha DP_{k_1}^{D-1}}^{1/\alpha DP_{k_1+1}^{D-1}} (
\alpha Ds)^{\beta -1}(P_{k-1}\cdots P_{k_1+1})^{D-1}s^{-\beta}
\,ds
\nonumber
\\[-8pt]
\\[-8pt]
\nonumber
& &\qquad\le(D-1) (\alpha D)^{\beta-1}(P_{k-1}\cdots
P_{k_1+1})^{D-1}\log(P_{k_1}/P_{k_1+1}).
\end{eqnarray}

By comparing terms involving $P_i$ and ignoring the other coefficients,
one can check that the largest
bound is given in (\ref{eq3145}). Since the logarithm term there is
dominated by $P_{k_1+1}^{-\delta(D-1)}$,
for given $\delta> 0$ and small enough $P_{k_1+1}$, it follows that~(\ref{eq364}) holds for $P_k$, for appropriate
$\fff$ and all $k$.
\end{pf*}

\section{\texorpdfstring{The case where $\beta\in(1,D/(D-1))$}
{The case where beta in (1,D/(D-1))}}\label{sec4}

In this section, we demonstrate Theorem~\ref{thm131}. We do this by
demonstrating the lower
and upper bounds needed for the theorem in Propositions~\ref{prop411}
and~\ref{prop441}.
Here, we set
\[
\nu_{\beta} = (\beta-1)/\bigl[1 - (D-1)(\beta- 1)\bigr].
\]

%
\begin{prop}
\label{prop411}
Consider a family of JSQ networks, with given $D\ge2$ and
$N=D,D+1,\ldots,$ where the $N$th network
has Poisson rate-$\alpha N$ input, with $\alpha< 1$, and where service
at each queue is FIFO, with
distribution $F(\cdot)$ having mean $1$. Assume that (\ref{eq112})
holds and that
%
%
\begin{equation}
\label{eq412} %
\bar{F}(s) \ge s^{-\beta} \qquad\mbox{for } s\ge
s_0, 
\end{equation}
with $\beta\in(1, D/(D-1))$ and some $s_0 \ge1$. Then, for
appropriate $\aaa> 0$
and all $k$,
%
%
\begin{equation}
\label{eq413} %
P_k \ge\aaa k^{-\nu_{\beta}}. 
\end{equation}
\end{prop}

%
\begin{prop}
\label{prop441}
Consider a family of JSQ networks, with given $D\ge2$ and
$N=D,D+1,\ldots,$ where the $N$th network
has Poisson rate-$\alpha N$ input, with $\alpha< 1$, and where service
at each queue is FIFO, with
distribution $F(\cdot)$ having mean $1$. Assume that (\ref{eq112})
holds and that
%
%
\begin{equation}
\label{eq442} %
\bar{F}(s) \le
s^{-\beta}\qquad \mbox{for } s\ge s_0, 
\end{equation}
with $\beta\in(1, D/(D-1))$ and some $s_0 \ge1$. Then, for each
$\delta> 0$,
appropriate $\doo> 0$, and all $k$,
%
%
\begin{equation}
\label{eq443} %
P_k \le\doo k^{-(1-\delta)\nu_{\beta}}. 
\end{equation}
\end{prop}
Theorem~\ref{thm131} follows immediately from Propositions \ref
{prop411} and~\ref{prop441}
upon letting $\delta\searrow0$ in (\ref{eq443}).

As in Section~\ref{sec3}, the demonstration of the lower bound is much quicker
than that of the upper bound. We
first demonstrate the lower bound, Proposition~\ref{prop411}, and
then, in the remainder of the section,
derive the upper bound, Proposition~\ref{prop441}.

\subsection*{Demonstration of Proposition \protect\ref{prop411}}
As in Section~\ref{sec3}, when we considered the case where $\beta> D/(D-1)$,
for the lower bound, it suffices to
construct a path along which $Z^{\mathcal{H}}(t)$ increases from $0$
to $k$ within the first cycle. As
before, we allocate the same amount of time for each of the first $k$
arrivals, which are also required
to occur before the first departure.

\begin{pf*}{Proof of Proposition~\ref{prop411}}
Consider the cavity process $X^{\mathcal{H}}(\cdot)$ with
$X^{\mathcal{H}}(0) = 0$.
We obtain a lower bound on the expected amount of time over which
$Z^{\mathcal{H}}(t) \ge k$
before $X^{\mathcal{H}}(\cdot)$ returns to $0$, assuming that $k\ge s_0$.

We consider the event where the first service time is at least $s_1 =
4k/\break(\alpha P_k^{D-1})$ and
the first $k$ arrivals occur by time $s_1/2$. We note that the
probability of the latter event occurring
is greater than the probability of at least $k$ events occurring by
time $s_1/2$ for a
rate-$\alpha P_k^{D-1}$ Poisson process, which, by a simple large
deviations estimate, is at least
\[
1 - e^{\bbb k} \ge1/2
\]
for large enough $k$ and an appropriate constant $\bbb$. Together with
(\ref{eq412}), this implies
that the expected amount of time in $[s_1/2,s_1]$, during which
$Z^{\mathcal{H}}(t) \ge k$ and
before $X^{\mathcal{H}}(\cdot)$ has returned to $0$, is at least
%
%
\begin{equation}
\label{eq414} %
\frac{1}{2} \cdot\frac{s_1}{2} \cdot
\bar{F}(s_1) \ge\frac{1}{4}\bigl(4k/\bigl(\alpha
P_k^{D-1}\bigr)\bigr)^{-(\beta- 1)}. %
\end{equation}

Inequality (\ref{eq414}) implies that
\[
P_k \ge\frac{\alpha}{16m_0}k^{-(\beta-1)} P_k^{(D-1)(\beta- 1)},
\]
where $m_0$ is the mean return time to $0$.
Solving for $P_k$, it follows from this that, for large $k$,
\[
P_k \ge\frac{\alpha}{16m_0} k^{-\nu_{\beta}},
\]
which implies (\ref{eq413}) for all $k$.
\end{pf*}

\subsection*{Demonstration of Proposition \protect\ref{prop441}}
The demonstration of the upper bound (\ref{eq443}) for Theorem \ref
{thm131} is considerably
more involved than is the lower bound. The basic idea is to consider
two cases, depending on
whether or not there is a service time~$s$ with $s > s_1$, for
preassigned $s_1 \ge1$, before
a state $x$ with $z=k$ is reached in the first cycle,
and to obtain upper bounds for each case.
The two bounds are given in Propositions~\ref{prop451} and \ref
{prop461}, which are then combined
in Corollary~\ref{cor471}. Employing Corollary~\ref{cor471}, the
proof of Proposition~\ref{prop441}
provides an iteration scheme where
a sequence of values $s_1(n), n=0,1,2,\ldots,$
for $s_1$ are given that provide successively better upper bounds for
$P_k$, and
that yield (\ref{eq443}) in the limit. The demonstration of
Proposition~\ref{prop461} involves
the construction of a supermartingale, whose details are postponed
until the end of the section.\vadjust{\goodbreak}

Let $\tau_k$, for given $k\in\mathbb{Z}_+$, denote the first time
$t$ in the first cycle at which
$Z^{\mathcal{H}}(t)=k$.
For Propositions~\ref{prop451} and~\ref{prop461}, we denote by
$B_{s_1,k}$ the set of realizations on which
some service time that is strictly greater than $s_1$, with $s_1 \ge
1$, occurs up to and including
the service time interval that contains $\tau_k$. Proposition~\ref{prop451}
considers the case where $B_{s_1,k}$ occurs; the demonstration of
the proposition is quick, using Lemma~\ref{lem391}. As in Sections~\ref{sec2}
and~\ref{sec3}, we
denote by $V_k$ the occupation time at states~$x$, with $z\ge k$.
%
%
\begin{prop}
\label{prop451}
Consider a family of JSQ networks with the same assumptions holding as
in Proposition~\ref{prop441}.
Then, for appropriate $\iii$ and all~$k$,
%
%
\begin{equation}
\label{eq452} %
E[V_k; B_{s_1,k}] \le\iii
s_1^{-\beta}(k + s_1). %
\end{equation}
\end{prop}
\begin{pf}
We apply Lemma~\ref{lem391} at the beginning of the first service time that
is greater than $s_1$. Since there are less than $k$ jobs
under $B_{s_1,k}$ then, it follows that,
for appropriate $\jjj$ and large enough $k$,
%
%
\begin{eqnarray}
\label{eq453} %
E[V_k; B_{s_1,k}] &\le&3
\bigl(P(B_{s_1,k})/\bar{F}(s_1)\bigr)\int
_{s_1}^{\infty}(k+s) F(ds)
\nonumber
\\[-8pt]
\\[-8pt]
\nonumber
&\le&\jjj\int_{s_1}^{\infty}(k+s) F(ds). %
\end{eqnarray}
For the latter inequality, note that there are only a finite expected
number of service times in
the first cycle, and that, by Wald's equation, the expected number of
such times that are at most $s$,
for given $s \ge0$, is proportional to $F(s)$. Since $k+s$ is
increasing in $s$, integration by
parts together with (\ref{eq442}) implies that the last quantity in
(\ref{eq453}) is at most
$\iii s_1^{-\beta}(k+s_1)$, for appropriate $\iii$.
\end{pf}

In order to consider the behavior of $X^{\mathcal{H}}(\cdot)$ on
$B_{s_1,k}^c$, we find it convenient to
employ the service time distribution $F^{s_1}(\cdot)$ that is given by
%
%
\begin{eqnarray}
\label{eq454} %
F^{s_1}(s) &=& F(s)\qquad \mbox{for } s <
s_1,
\nonumber
\\[-8pt]
\\[-8pt]
\nonumber
&=& 1\qquad \mbox{for } s\ge s_1. %
\end{eqnarray}
We define $X_{s_1}^{\mathcal{H}}(\cdot)$ analogously to $X^{\mathcal
{H}}(\cdot)$, but where the
service time distribution of the process is $F^{s_1}(\cdot)$ up to and
including the service time
interval containing $\tau_k$, and is given by $F(\cdot)$ afterwards;
$Z_{s_1}^{\mathcal{H}}(\cdot)$ and $S_{s_1}^{\mathcal{H}}(\cdot)$
are defined analogously.
One has
%
%
\begin{equation}
\label{eq455} %
E\bigl[V_k; B_{s_1,k}^c
\bigr] \le E\bigl[V_k^{s_1}\bigr], %
\end{equation}
where $V_k^{s_1}$ is the occupation time at states $x$ with $z\ge k$ for
$X_{s_1}^{\mathcal{H}}(\cdot)$. Note that the mean of $F^{s_1}(\cdot
)$ is at most $1$.

In contrast to Proposition~\ref{prop451}, Proposition~\ref{prop461} requires us
to restrict our
choice of $s_1$ in terms of $k$. For this, we set $k_1 = \lfloor k/3
\rfloor$ and\vadjust{\goodbreak} introduce the abbreviation
%
%
\begin{equation}
\label{eq464} %
p = p_{k_1} = \alpha D P_{k_1}^{D-1}.
\end{equation}
The required restriction on $s_1$ is that
%
%
\begin{equation}
\label{eq467} %
s_1 \le k^{1-\eta}/p, %
\end{equation}
where $\eta\in(0,1/2)$.
In the proof of Proposition~\ref{prop441}, we will introduce an
iterative scheme that involves explicit choices
of $s_1$ based on our knowledge of $P_{k_1}$ at each step.

Proposition~\ref{prop461} gives us the following upper bound for
$E[V_k; B_{s_1,k}^c]$.
%
%
\begin{prop}
\label{prop461}
Consider a family of JSQ networks with the same assumptions holding as
in Proposition~\ref{prop441}.
Suppose that $\delta>0$ and $\eta\in(0,1/2)$ are given, and that
$s_1$ satisfies (\ref{eq467}).
Then, for appropriate $\kkk$ and all $k$,
%
%
\begin{equation}
\label{eq462} %
E\bigl[V_k; B_{s_1,k}^c
\bigr] \le\kkk(k + s_1) \exp\bigl\{-\delta k^{\eta}\bigr
\}. %
\end{equation}
\end{prop}

\begin{figure}

\includegraphics{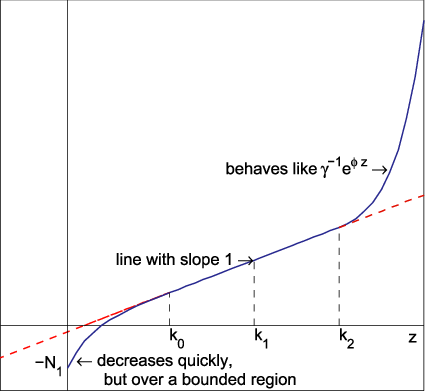}

\caption{Graph of $f(z)$.}\label{figurenew2f}
\end{figure}

The demonstration of Proposition~\ref{prop461} depends on an appropriate
supermartingale. In order to construct the supermartingale, we employ
the following notation.
We fix $k_0 \in\mathbb{Z}_+$, which will not depend on $k$ as $k$
increases, and set
$k_2 = 2k_1$, where~$k_1$ is as defined earlier. We set
%
%
\begin{eqnarray}
\label{eq463} %
f(z) &=& (z\wedge k_2) - N_1
\exp\bigl\{-\theta(z\wedge k_0)\bigr\}
\nonumber
\\[-8pt]
\\[-8pt]
\nonumber
&&{} + \gamma^{-1}\exp
\bigl\{\phi(z\vee k_2)\bigr\} - \gamma^{-1}\exp\{\phi
k_2\}, %
\end{eqnarray}
where $N_1,\theta> 0$, $\phi= \delta k^{\eta- 1}$ and $\gamma=
\phi e^{\phi k_2}$, and where $\delta>0$ and $\eta\in(0,1/2)$
are as in Proposition~\ref{prop461}; the function $f(\cdot)$ is
sketched in Figure~\ref{figurenew2f}.
The terms $P_k$ will continue to refer to the probabilities defined at
the beginning of the paper with
respect to the cavity process with the original service distribution
$F(\cdot)$ [not $F^{s_1}(\cdot)$].

%

We let $H(n)$, with $n\ge1$, denote the number of jobs for the process
$X_{s_1}^{\mathcal{H}}(\cdot)$, with $X_{s_1}^{\mathcal{H}}(0)=0$,
at the time when the $n$th job has just departed; we set $H(0)=1$, and
we let $\rho$
denote the first time $n$ at which either $H(n)=0$ or $H(n)\ge k-1$.
Using this notation, we define the analog of $M(\cdot)$ in
(\ref{eq3101}),
%
%
\begin{equation}
\label{eq465} %
M(n) = f\bigl(H(n\wedge\rho)\bigr). %
\end{equation}
Note that, unlike for $M(\cdot)$ in (\ref{eq3101}), $M(\cdot)$ here
depends strongly on the choice of~$k$.
Also, unlike $M(\cdot)$ in (\ref{eq3101}), it was not necessary to
wait until the first departure in defining $H(0)$,
since $X_{s_1}^{\mathcal{H}}(0)=0$, and hence there is no initial
residual service time; in both cases,
$H(1)-H(0)$ is the change in the number of jobs during the service time
of the first job that begins service
when $t > 0$.
%
%
\begin{prop}
\label{prop466}
Consider a family of JSQ networks with the same assumptions holding as
in Proposition~\ref{prop441}.
Suppose that $\delta>0$ and $\eta\in(0,1/2)$ are given, and that
$M(\cdot)$ is defined as above.\vadjust{\goodbreak}
Also, assume that $s_1$ satisfies (\ref{eq467}).
Then, for large enough $k$, $M(\cdot)$ is a supermartingale, with respect
to the filtration $\mathcal{G}_n = \sigma(H(0),\ldots,H(n))$, for
small enough $\delta>0$, and appropriate $\theta,N_1 > 0$, with
$\delta$, $\theta$, and $N_1$ not depending on $k$.
\end{prop}

The demonstration of Proposition~\ref{prop466} will be given at the
end of the section.
Employing Proposition~\ref{prop466}, we now demonstrate Proposition
\ref{prop461}.
\begin{pf*}{Proof of Proposition~\ref{prop461}}
We suppose that the terms $\delta$, $\theta$ and $N_1$ are chosen so
that, for large enough
$k$, $M(\cdot)$ is a supermartingale.
Set $\sigma_{L} = \min\{n\dvtx M(n) \ge L\}$, for given $L >0$, which
will depend on $k$.
Since $M(\cdot)$ is bounded below by $-N_1$ and $M(0)\le1$, by the
optional sampling theorem,
%
%
\begin{equation}
\label{eq468} %
P(\sigma_{L} < \infty) \le
\frac{1}{L}(1 + N_1). %
\end{equation}

On the other hand, denoting by $n_k$ the service interval during which
$Z_{s_1}^{\mathcal{H}}(t) = k$ first occurs and by $T_{n_k}$ the end
of that interval,
$H(n_k) = Z_{s_1}^{\mathcal{H}}(T_{n_k}) \ge k-1$. Substituting this
into (\ref{eq463})--(\ref{eq465})
and recalling that $\phi= \delta k^{\eta-1}$, one obtains
\[
M(n_k) \ge-N_1 + \gamma^{-1}\exp\bigl\{
\phi(k-1)\bigr\} - \gamma^{-1}\exp\{2\phi k/3\} \ge\exp\bigl\{
\delta
k^{\eta}\bigr\}/2\gamma
\]
for large $k$.
Let $\tau_k^{s_1}$ denote the first time $t$, during the first cycle,
at which $Z_{s_1}^{\mathcal{H}}(t) = k$. Plugging
$L = \exp\{\delta k^{\eta}\}/2\gamma$ into (\ref{eq468}),
substituting in for $\gamma$ and recalling that $k_2 =
2\lfloor k/3 \rfloor$, it follows that, for large $k$,
%
%
\begin{eqnarray}
\label{eq4610} %
P\bigl(\tau_k^{s_1} < \infty
\bigr) \le P(\sigma_L < \infty) &\le&\exp\bigl\{-\delta
k^{\eta}\bigr\} \cdot\exp\bigl\{2\delta k^{\eta}/3\bigr\}
\nonumber
\\[-8pt]
\\[-8pt]
\nonumber
&= &\exp\bigl\{-\delta k^{\eta}/3\bigr\}. %
\end{eqnarray}

Lemma~\ref{lem391} applied to $F(\cdot)$, which is the service
distribution of new service times
after $\tau_k^{s_1}$, provides the upper bound
\[
E\bigl[V_k^{s_1} | \mathcal{F}_{\tau_k^{s_1}}\bigr] \le
2(k+s+N_0),
\]
given that $S_{s_1}^{\mathcal{H}}(\tau_k^{s_1})=s$. Since the
residual service time for $X_{s_1}^{\mathcal{H}}(t)$
is at most $s_1$ for $t\le\tau_k^{s_1}$, it therefore follows from
(\ref{eq4610}) that, for large $k$,
%
%
\begin{equation}
\label{eq4612} %
E\bigl[V_k^{s_1}\bigr] \le
3(k+s_1)\exp\bigl\{-\delta k^{\eta}/3\bigr\}. %
\end{equation}
The inequality in (\ref{eq462}) follows upon applying (\ref{eq455})
to (\ref{eq4612}) and
substituting in a smaller choice of $\eta$.
\end{pf*}

We combine the upper bounds given in Propositions~\ref{prop451} and
\ref{prop461} for $E[V_k; B_{s_1}]$
and $E[V_k; B_{s_1}^c]$ to obtain the following upper bound on
$E[V_k]$. Since we will always assume
$s_1 \le k^{\nu_{\beta}+1}$ in our application of the corollary,
this allows us to omit the exponential term inherited from (\ref{eq462}).
%
%
\begin{cor}
\label{cor471}
Consider a family of JSQ networks with the same assumptions holding as
in Proposition~\ref{prop441}.
Fix $\eta\in(0,1)$ and assume that
%
%
\begin{equation}
\label{eq472} %
s_1 \le\bigl[(\alpha
D)^{-1}k^{1-\eta} P_{k_1}^{1-D}\bigr] \wedge
k^{N} 
\end{equation}
for some $N>0$. Then, for appropriate $\nnn$ and all $k$,
%
%
\begin{equation}
\label{eq473} %
E[V_k] \le\nnn s_1^{-\beta}(k+s_1).
\end{equation}
\end{cor}
\begin{pf}
It follows from Propositions~\ref{prop451} and~\ref{prop461} that
\[
E[V_k] \le\iii s_1^{-\beta}(k+s_1)
+ \kkk(k+s_1)\exp\bigl\{-\delta k^{\eta}\bigr\}
\]
for appropriate $\iii$ and $\kkk$. The assumption $s_1 \le k^{N}$
allows us
to absorb the second term into the first.
\end{pf}

The following elementary lemma will be employed in the proof of
Proposition~\ref{prop441}.
%
%
\begin{lem}
\label{lem481}
Suppose that R(n) satisfies
%
%
\begin{equation}
\label{eq482} %
R(n) = aR(n-1) + b\qquad \mbox{for } n\ge1, %
\end{equation}
with $R(0) = c$, for $a\in(0,1)$ and $b,c \in\mathbb{R}$. Then,
%
%
\begin{equation}
\label{eq483} %
\lim_{n\rightarrow\infty}R(n) = b/(1-a). %
\end{equation}
If $R(0) < b/(1-a)$, then the sequence $R(n)$ is increasing, and if
$R(0) > b/(1-a)$, then
the sequence is decreasing.
\end{lem}
\begin{pf}
Setting $\tilde{R}(n) = R(n) - b/(1-a)$, it follows from (\ref
{eq482}) that
%
%
\begin{equation}
\label{eq484} \tilde{R}(n) = a\tilde{R}(n-1) \qquad\mbox{for } n\ge1,
\end{equation}
with $\tilde{R}(0) = c - b/(1-a)$. All of the claims follow by
iterating (\ref{eq484}).
\end{pf}
%
%
%
%
%
%
%

We will employ the lemma in the following multiplicative format.
%
%
\begin{cor}
\label{cor487}
Suppose that $Q_k(n)$ satisfies
%
%
\begin{equation}
\label{eq488} %
Q_k(n) = \bigl(k^{-(1-2\eta)}Q_k(n-1)^{D-1}
\bigr)^{\beta- 1}\qquad \mbox{for } n\ge1, %
\end{equation}
with $Q_k(0) = k^{1- \beta+ 2\eta\beta}$,\hspace*{2pt} for $(D-1)(\beta-1) \in
(0,1)$ and $\eta\in(0,1/2)$. Then,
$Q_k(n)$ satisfies $Q_k(n) = k^{-R(n)}$, where the sequence $R(n)$ is
increasing in $n$ and
%
%
\begin{equation}
\label{eq489} %
\lim_{n\rightarrow\infty}R(n) = (1-2\eta)
\nu_{\beta}, %
\end{equation}
with $\nu_{\beta} = (\beta-1)/[1-(D-1)(\beta-1)]$.
\end{cor}
\begin{pf}
The limit in (\ref{eq489}) follows from (\ref{eq483}) upon setting
$a = (D-1)(\beta-1)$, $b = (1-2\eta)(\beta-1)$ and $c = \beta- 1 -
2\eta\beta$. The sequence $R(n)$
is increasing since $R(0) < (1-2\eta)\nu_{\beta}$.
\end{pf}

We now employ Corollaries~\ref{cor471} and~\ref{cor487} to
demonstrate Proposition~\ref{prop441}.
\begin{pf*}{Proof of Proposition~\ref{prop441}}
For given $k$ and $\eta\in(0,1/2)$, we define $Q_k(n)$ as in
Corollary~\ref{cor487} and set
%
%
\begin{eqnarray}
\label{eq491} %
s_1(n) &=& (\alpha D)^{-1}k^{1-\eta}\qquad
\mbox{for } n=0,
\nonumber
\\[-8pt]
\\[-8pt]
\nonumber
&=& \bigl(\alpha D Q_{k_1}(n-1)^{D-1} \bigr)^{-1}k^{1-\eta}
\qquad\mbox{for } n\ge1, %
\end{eqnarray}
where $k_1 = \lfloor k/3 \rfloor$. Using $s_1(n)$, we will inductively
show that,
for large $k$ (depending on $\eta$),
%
%
\begin{equation}
\label{eq492} %
P_k \le Q_k(n) \qquad\mbox{for all }
n\ge0. %
\end{equation}
Letting $n\rightarrow\infty$, it therefore follows from the corollary that
%
%
\begin{equation}
\label{eq493} %
P_k \le k^{-(1-2\eta)\nu_{\beta}}. %
\end{equation}
This implies (\ref{eq443}) in Proposition~\ref{prop441}, with $\delta
< 2\eta$.

To show (\ref{eq492}) holds for $n=0$, we note that $s_1(0)$ satisfies
(\ref{eq472}).
Therefore, by (\ref{eq233}) and Corollary~\ref{cor471}, for large $k$,
%
%
\begin{equation}
\label{eq494} %
P_k \le2\nnn(m_0)^{-1}s_1(0)^{-\beta}k
\le k^{-(\beta-1)+2\eta\beta} = Q_k(0), %
\end{equation}
where the constants in the second expression are absorbed in the third
expression by using
the $2\eta$ term.
Note that, in this application of (\ref{eq473}), $s_1(0)\le k$. In the
application of
(\ref{eq473}) given next, $s_1(n)\ge k$ for all $n\ge1$.

Suppose that (\ref{eq492}) holds with $n-1$ in place of $n$. Choosing
$s_1(n)$ as in
(\ref{eq491}) and employing the lower bound for $Q_k(n)$ given in
(\ref{eq489}), one can
check that $s_1(n)$ satisfies (\ref{eq472}), with $N=\nu_{\beta}+1$.
Also note that, by Corollary~\ref{cor487},
\[
Q_{k_1}(n)\le3^{\nu_{\beta}}Q_{k}(n)
\]
for large $k$ and all $n$.
Applying (\ref{eq233}) and Corollary~\ref{cor471} again, we therefore
obtain that, for large $k$,
%
%
\begin{eqnarray}
\label{eq496} %
P_k \le2\nnn(m_0)^{-1}s_1(n)^{-(\beta-1)}
&\le&\bigl(k^{-(1-2\eta)}Q_k(n-1)^{D-1}
\bigr)^{\beta-1}
\nonumber
\\[-8pt]
\\[-8pt]
\nonumber
& =& Q_k(n). %
\end{eqnarray}
This demonstrates (\ref{eq492}).
\end{pf*}
%
%
%
%
%
%

In order to complete the demonstration of Proposition~\ref{prop441},
we need to prove
Proposition~\ref{prop466}, which asserts that $M(\cdot)$, given by
(\ref{eq465}), is
a supermartingale.

\begin{pf*}{Proof of Proposition~\ref{prop466}}
We need to show the supermartingale inequality (\ref{eq3102}) for
$H(n)\in(0,k-1)$. We do this separately over the intervals $(0,k_1]$
and $(k_1,k-1)$.
The basic idea for the first interval will be to show
that, on $(0,k_1]$, (\ref{eq3102}) will be satisfied for the same
reasons as was $M(\cdot)$, for
$M(\cdot)$ given by (\ref{eq3101}), the point being that, since $k_2
- k_1 = \lfloor k/3 \rfloor$
is large, the role played by the additional terms $\gamma^{-1}\exp\{
\phi(z\vee k_2)\}-\gamma^{-1}\exp\{\phi k_2\}$ in
(\ref{eq463}) is negligible. On the second interval $(k_1,k-1)$, the
strong negative drift of
$Z_{s_1}^{\mathcal{H}}(\cdot)$ will be enough to compensate for both
the $z\wedge k_2$ and
$\gamma^{-1}\exp\{\phi(z\vee k_2)\}-\gamma^{-1}\exp\{\phi k_2\}$
terms. We do the latter interval first.

We claim that for large $k$ and $H(n)\ge k_1$,
%
%
\begin{equation}
\label{eq4101} %
E\bigl[\exp\bigl\{\phi H(n+1)\bigr\}|
\mathcal{G}_n\bigr] \le E\bigl[\exp\bigl\{\phi H(n)\bigr\}\bigr]. %
\end{equation}
We first note that,
because of (\ref{eq464}), for $H(n)\ge k_1$, the number of arrivals
over the $(n+1)$st service
interval is dominated by a mixture of Poisson rate-$ps$ random
variables, with $s$ being distributed
according to $F^{s_1}(\cdot)$. Therefore,
\[
E\bigl[\exp\bigl\{\phi\bigl(H(n+1)-H(n)\bigr)\bigr\}|
\mathcal{G}_n\bigr] \le e^{-\phi}\int_{0}^{s_1}
\exp\bigl\{ps\bigl(e^{\phi} -1\bigr)\bigr\}F^{s_1}(ds).
\]
Since the integrand is convex and the mean of $F^{s_1}(\cdot)$ is at
most $1$, the right-hand side is at most
%
%
\begin{equation}
\label{eq4102} %
e^{-\phi} \biggl[ \biggl(1-\frac{1}{s_1}
\biggr) + \frac{1}{s_1}\exp\bigl\{ps_1 \bigl(e^{\phi}-1
\bigr)\bigr\} \biggr]. %
\end{equation}
On account of the definitions of $\phi$ and $p$ given between (\ref
{eq464}) and (\ref{eq465}), both
$\phi$ and $ps_1 \phi$ are at most\vadjust{\goodbreak} $\delta$. Using $e^z \sim1+z$
for $z$ close to $0$,
one can therefore check that, for given $\varepsilon> 0$ and small
enough $\delta> 0$, (\ref{eq4102})
is at most
\[
1 + \phi\bigl[(1+\varepsilon)p - (1-\varepsilon)\bigr]. %
\]
For $p \le(1-\varepsilon)/(1+\varepsilon)$, the above quantity is at
most $1$, which holds here
since $p\rightarrow0$ as $k\rightarrow\infty$. This implies (\ref{eq4101}).

For $H(n) > k_2$, it is easy to see that (\ref{eq3102}) follows from
(\ref{eq4101}),
since
%
%
\begin{eqnarray}
\label{eq4104} %
f(z) - \gamma^{-1}e^{\phi z} &=& b\qquad
\mbox{for } z\ge k_2,
\nonumber
\\[-8pt]
\\[-8pt]
\nonumber
&\le& b\qquad \mbox{for } z< k_2, %
\end{eqnarray}
where $b \stackrel{\mathrm{def}}{=} f(k_2) - \gamma^{-1}e^{\phi
k_2}$. For $H(n) \in(k_1,k_2]$,
(\ref{eq3102}) follows from (\ref{eq4101}) with a bit more work. In
place of (\ref{eq4104}),
one uses
%
%
\begin{equation}
\label{eq4105} %
g(z)\stackrel{\mathrm{def}} {=} f(z) -
\gamma^{\prime}e^{\phi z} \le f\bigl(H(n)\bigr) -
\gamma^{\prime} e^{\phi H(n)} %
\end{equation}
for all $z$, where $\gamma^{\prime} \stackrel{\mathrm{def}}{=}
(\phi e^{\phi H(n)})^{-1} =
\gamma^{-1}e^{\phi(k_2 - H(n))}$. To check (\ref{eq4105}), note that
equality holds for
$z=H(n)$; we claim that the maximum of $g(\cdot)$ is taken there. One
has $g^{\prime}(H(n)) = 0$
because of our definition of
$\gamma^{\prime}$; $g^{\prime}(z) \ge0$ for $z\le H(n)$ and
$g^{\prime}(z) \le0$ for $z\in[H(n),k_2)$
because of the concavity of $g(\cdot)$ there; and since $\gamma
^{\prime} \ge1$, for $z> k_2$,
it is easy to see that $g^{\prime}(z) \le0$ there. This shows (\ref
{eq4105}) and hence (\ref{eq3102})
for $H(n) \in(k_1,k_2]$ as well.

We still need to show (\ref{eq3102}) for $H(n)\in(0,k_1]$. For this,
we compare $M(\cdot)$ with
$\tilde{M}(\cdot)$, where
\[
\tilde{f}(z) = z + n/2 - N_1 \exp\bigl\{-\theta(z\wedge
k_0)\bigr\}
\]
and
\[
\tilde{M}(n) = \tilde{f}\bigl(H(n\wedge\rho)\bigr).
\]
Set $R(n) = M(n) - \tilde{M}(n)$. For $H(n)\in(0,k_1]$, one has
%
%
\begin{eqnarray}
\label{eq4108} %
R(n+1)-R(n) + 1/2 &=& 0 \qquad\mbox{for } H(n+1) \le
k_2,
\nonumber
\\[-8pt]
\\[-8pt]
\nonumber
&\le&\gamma^{-1}e^{\phi H(n+1)} \qquad\mbox{for } H(n+1) > k_2.
\end{eqnarray}
Since $\tilde{M}(\cdot)$ is the supermartingale in (\ref{eq3101}),
except with a different
initial state, $\tilde{M}(\cdot)$ satisfies (\ref{eq3102}) if
$\theta$ and $N_1$ are chosen
as in (\ref{eq3101}).
In a moment, we will show that
%
\begin{equation}
\label{eq4109} %
E\bigl[e^{\phi H(n+1)}\mathbf{1}\bigl
\{H(n+1)>k_2\bigr\}| \mathcal{G}_n\bigr] \le\gamma/2
\end{equation}
for $H(n) \le k_1$ and large $k$. Using (\ref{eq4108}) and (\ref
{eq4109}), (\ref{eq3102})
therefore also follows for $M(\cdot)$ for $H(n) \le k_1$.

It suffices to show (\ref{eq4109}) for $H(n) = k_1$. To do this,
we need to control the right tail of $H(n+1)$. The number of
arrivals over the $(n+1)$st service interval for the cavity process is
dominated by a mixture of
Poisson mean-$ps_1$ random variables, with the mixture distributed
according to $F^{s_1}$. This
mixture is in turn dominated by a Poisson mean-$s_1$ random variable.
Therefore, the left-hand side
of (\ref{eq4109}) is at most
%
%
\begin{equation}
\label{eq41010} %
\sum_{k^{\prime}=k_2}^{\infty}
\bigl[e^{-ps_1}(ps_1)^{k^{\prime}-k_1}/\bigl(k^{\prime}-k_1
\bigr)! \bigr]e^{\phi k^{\prime}}. %
\end{equation}
Setting $\ell= k^{\prime} - k_2$, one has
\[
\bigl(k^{\prime}-k_1\bigr)! \ge\ell! (k_2-k_1)!
\ge\ell! \bigl((k_2-k_1)/e\bigr)^{k_2-k_1},
\]
where the last inequality follows from Stirling's formula.
Substituting $\ell$ into (\ref{eq41010}), applying this bound, and employing
$\exp\{e^{\phi}ps_1\} = \sum_{\ell= 0}^{\infty}(e^{\phi
}ps_1)^{\ell}/\ell!$,
it follows that (\ref{eq41010}) is at most
%
%
\begin{equation}
\label{eq41011} 
\biggl(\frac{eps_1}{k_2-k_1}
\biggr)^{k_2-k_1} \exp\bigl\{ps_1\bigl(e^{\phi} -1
\bigr) + \phi k_2 \bigr\} \le\mmm k^{-\eta k/3} e^{4\phi k}
\end{equation}
for appropriate $\mmm$, where the inequality employs (\ref{eq467})
and $e^{\phi} -1 \le2\phi$, for small $\phi$.
As $k\rightarrow\infty$, the right-hand side of (\ref{eq41011}) goes
to $0$. It follows that
the left-hand side of (\ref{eq4109}), with $H(n) = k_1$, goes to $0$
as $k\rightarrow\infty$.
This implies (\ref{eq4109}) holds for $H(n)\le k_1$ and large $k$,
which completes the proof of the proposition.
\end{pf*}

\section{\texorpdfstring{The case where $\beta= D/(D-1)$}
{The case where beta = D/(D-1)}}\label{sec5}

In this section, we demonstrate Theorem~\ref{thm141}.
We do this by demonstrating the lower and upper
bounds needed for the theorem, in Propositions~\ref{prop511} and \ref
{prop531}.
%
%
\begin{prop}
\label{prop511}
Consider a family of JSQ networks, with given $D\ge2$ and
$N=D,D+1,\ldots,$ where the $N$th network
has Poisson rate-$\alpha N$ input, with $\alpha< 1$, and where service
at each queue is FIFO, with
distribution $F(\cdot)$ having mean $1$. Assume that (\ref{eq112})
holds and that
%
%
\begin{equation}
\label{eq512} %
\bar{F}(s) \ge c_1 s^{-D/(D-1)}\qquad
\mbox{for } s\ge s_0, 
\end{equation}
for some $c_1 > 0$ and $s_0 \ge1$.
Then, for appropriate $\ccc> 0$ and $s_D(c_1) < \infty$,
%
%
\begin{equation}
\label{eq513} %
P_k \ge\ccc e^{-s_D(c_1) k} \qquad\mbox{for
all } k, 
\end{equation}
where
%
%
\begin{equation}
\label{eq514} %
s_D(c_1) \searrow0\qquad
\mbox{as } c_1\nearrow\infty. %
\end{equation}
\end{prop}
%

%
\begin{prop}
\label{prop531}
Consider a family of JSQ networks, with given $D\ge2$ and
$N=D,D+1,\ldots,$ where the $N$th network
has Poisson rate-$\alpha N$ input, with $\alpha< 1$, and where service
at each queue is FIFO, with
distribution $F(\cdot)$ having mean $1$. Assume that (\ref{eq112})
holds and that
%
%
\begin{equation}
\label{eq532} %
\bar{F}(s) \le c_2 s^{-D/(D-1)}\qquad
\mbox{for } s\ge s_0,\vadjust{\goodbreak} 
\end{equation}
for some $c_2 < \infty$ and $s_0 \ge1$.
Then, for appropriate $\ooo$ and $r_D(c_2) > 0$,
%
%
\begin{equation}
\label{eq533} %
P_k \le\ooo e^{-r_D(c_2) k}\qquad \mbox{for
all } k, 
\end{equation}
where
%
%
\begin{equation}
\label{eq534} %
r_D(c_2) \nearrow\infty\qquad
\mbox{as } c_2 \searrow0. %
\end{equation}
\end{prop}

Theorem~\ref{thm141} follows immediately from Propositions~\ref{prop511}
and~\ref{prop531}.

As in the previous two sections, the demonstration of the lower bound is
substantially quicker than that of the upper bound. We first demonstrate
the lower bound, Proposition~\ref{prop511} and then, in the remainder of
the section, derive the upper bound Proposition~\ref{prop531}.

\subsection*{Demonstration of Proposition \protect\ref{prop511}}
As in Sections~\ref{sec3} and~\ref{sec4}, where we considered the cases $\beta> D/(D-1)$ and
$\beta\in(1,D/(D-1))$, for the lower bound, it suffices to construct
a path
along which $Z^{\mathcal{H}}(t)$ increases from $0$ to $k$ within the first
cycle. In contrast to the previous two settings, we allocate geometrically
increasing amounts of time to the sequence of arrivals, up through the $k$th
arrival; as before, these arrivals are required to occur before the
time of
the first departure.

\begin{pf*}{Proof of Proposition~\ref{prop511}}
The argument is similar to that for Proposition~\ref{prop411} in that we
examine the cavity process $X^{\mathcal{H}}(\cdot)$ with
$X^{\mathcal{H}}(0)=0$, and obtain a lower bound on the expected amount
of time $V_k$ over which $Z^{\mathcal{H}}(t)\ge k$ before
$X^{\mathcal{H}}(\cdot)$ returns to $0$. Here, we argue by induction, and
assume that
%
%
\begin{equation}
\label{eq521} %
P_i \ge\ccc e^{-\aba i}\qquad \mbox{for }
i=0,\ldots, k-1, %
\end{equation}
for given $k$, where $\ccc\le[(\aba\vee1)s_0]^{-1}$, and $\aba>0$ will
be specified later.

We consider the following event $A$ that leads to a lower bound on
$P_k$ that
is compatible with (\ref{eq521}). We stipulate that the first service time
is at least
%
%
\begin{equation}
\label{eq522} s_1 \stackrel{\mathrm{def}} {=} \ddd
e^{\aba(D-1)k},
\end{equation}
where $\ddd= 4(\alpha\aba)^{-1}\ccc^{-(D-1)}$. Note that $\ddd\ge
s_0$. We also assume
that the interarrival time for the $(i+1)$st arrival at the queue,
$i=0,\ldots,k-1$,
is at most
%
%
\begin{equation}
\label{eq523} \alpha^{-1} \ccc^{-(D-1)}\exp\bigl\{
\tfrac{1}{2}\aba(D-1) (k+i)\bigr\}.
\end{equation}
A little estimation shows that the sum of the terms in (\ref{eq523}),
over $i=0,\ldots,k-1$, is bounded above by
%
%
\begin{eqnarray}
\label{eq524} %
&&\alpha^{-1} \ccc^{-(D-1)}\exp\bigl\{
\aba(D-1)k\bigr\}/\bigl(\exp\bigl\{\tfrac{1}{2}\aba(D-1)\bigr
\} - 1
\bigr)
\nonumber
\\[-8pt]
\\[-8pt]
\nonumber
&&\qquad \le(2/\alpha\aba) \ccc^{-(D-1)}\exp\bigl\{\aba(D-1)k\bigr\},
\end{eqnarray}
which is one-half of (\ref{eq522}).\vadjust{\goodbreak}

On account of the induction hypothesis in (\ref{eq521}), the
probability that the
$(i+1)$st arrival occurs within the interarrival time in (\ref{eq523})
is at least
\[
1 - \exp\bigl\{-e^{({1}/{2})\aba(D-1)(k-i)}\bigr\}. %
\]
%
%
%
%
So, the probability that the corresponding events for $i=0,\ldots,D-1$
all occur within
the allotted time is at least
\[
\prod_{i=1}^k \bigl(1 - \exp
\bigl\{-e^{({1}/{2})\aba(D-1)i} \bigr\} \bigr) \ge\psi(\aba),
\]
where $\psi(\aba) > 0$ for $\aba> 0$ and does not depend on $k$ or $D$,
with $\psi(\aba)\rightarrow1$ as $\aba\rightarrow\infty$; the
inequality requires
a little computation.

It follows from the previous two paragraphs that the event $A$, given
by the service time and interarrival
times restricted as in (\ref{eq522}) and (\ref{eq523}), has
probability at least
\[
\psi(\aba) \bar{F}\bigl(\ddd\exp\bigl\{\aba(D-1)k\bigr\}\bigr).
\]
On $A$, $Z^{\mathcal{H}}(t) \ge k$ over the interval $[s_1/2,s_1]$,
which has length
$\frac{1}{2}\ddd\exp\{\aba(D-1)k\}$. So,
%
\[
E[V_k] \ge\tfrac{1}{2}\ddd\psi(\aba) \exp\bigl\{
\aba(D-1)k\bigr\} \bar{F}\bigl(\ddd\exp\bigl\{\aba(D-1)k\bigr\}
\bigr).
\]
%
%
%
By substituting the bound in (\ref{eq512}) for $\bar{F}(s)$ and employing
$P_k = m_0^{-1}E[V_k]$, one obtains
\begin{eqnarray*}
P_k &\ge&\frac{1}{2m_0}\psi(\aba) c_1
\ddd\exp\bigl\{\aba(D-1)k\bigr\} \bigl(\ddd\exp\bigl\{\aba
(D-1)k\bigr\}
\bigr)^{-D/(D-1)}
\\
&=& \frac{1}{2m_0}\psi(\aba) c_1 (\ddd)^{-1/(D-1)}e^{-\aba k}
\ge\frac{c_1}{4m_0}\psi(\aba) (\alpha\aba)^{1/(D-1)} \ccc
e^{-\aba k}. %
\end{eqnarray*}

For given $c_1$ and large enough $\aba$, the last quantity
in the above display is at least $\ccc e^{-\aba k}$.
This implies the induction hypothesis in (\ref{eq521}) for $k$ and
this choice of $\aba$. Since
(\ref{eq521}) obviously holds for $i=0$, (\ref{eq513}) follows, with
$s_D(c_1) = \aba$. Similarly, for given $\aba$, one obtains the lower
bound $\ccc e^{-\aba k}$, if $c_1$ is chosen large enough, which implies
(\ref{eq514}). This completes the proof.
\end{pf*}

\subsection*{Demonstration of Proposition \protect\ref{prop531}}
The demonstration of the upper bound (\ref{eq533}) is substantially
more involved than is the lower
bound. The basic idea is similar to that employed for the upper bound
in Section~\ref{sec3}, where we
classified different paths for attaining $Z^{\mathcal{H}}(t)+k$, for
given $k$ and some $t$, in terms of
partitions $\pi$ given by (\ref{eq3111}). There, the probability of
the event associated with the
trivial partition dominated the probabilities for the other partitions.
Computing an upper bound
for the probability for the trivial partition and multiplying by the
upper bound $2^{\beta}$ for
the total number of partitions gave us our desired upper bounds on $P_k$.\vadjust{\goodbreak}

The details of our setup here will be different. The partitions we
consider will be
defined somewhat differently, and we will need to be more careful in
summing up probabilities---we will compute the probability of the event associated with the
trivial partition separately,
and will then sum up the probabilities for the other partitions, which
will be negligible in comparison.
We will also require an upper bound on $P_k$ from
Proposition~\ref{prop441}, at the beginning of our argument.
On the other hand, the computations of these upper bounds will be
substantially easier here than
the corresponding bounds were in Section~\ref{sec3}. The key difference is that
here the probabilities $P_k$ will decrease sufficiently
slowly in $k$ so that, for our estimates, not too much will be lost if
we consider $P_i$ to be
approximately the same for $i=k_1,\ldots,k-1$,
which will simplify our computations.

%
%
%
%
In order to show (\ref{eq533}) and (\ref{eq534}) of Proposition \ref
{prop531},
we will argue by induction, assuming that, for preassigned $\aca,\ooo
>0$ and $k_0,h_T \in\mathbb{Z}_+$,
%
%
\begin{equation}
\label{eq541} %
P_i \le\ooo e^{-\aca i}\qquad \mbox{for }
i=k_0,\ldots,k-1, %
\end{equation}
for given $k$ with $k\ge k_0 + h_T$. For appropriate choices of these
preassigned values, we will show
that the inequality in (\ref{eq541}) holds with $i=k$. We set
%
%
\begin{equation}
\label{eq541a} %
h_T = \bigl\lceil700 D^2
c_2 \bigr\rceil^{D-1} \vee6 %
\end{equation}
and
%
%
\begin{equation}
\label{eq541b} %
\aca= (h_T)^{-1} \vee
\tfrac{1}{6}\log\bigl(\bigl(220D^2 c_2
\bigr)^{-1} \bigr), %
\end{equation}
where $c_2$ is as in (\ref{eq532}). These particular choices of $h_T$
and $\aca$ are not
needed for most of the argument, and will only be inserted at the very end.

In order to specify
$\ooo$ and $k_0$, we note that, since (\ref{eq442}) is satisfied for
every $\beta< D/(D-1)$
because of (\ref{eq532}) and since $\nu_{\beta}\nearrow\infty$ as
$\beta\nearrow D/(D-1)$,
it follows from Proposition~\ref{prop441} that, for any $N$,
$\lim_{k\rightarrow\infty}k^N P_k = 0$.
Here, we set $N=h_T +1$. We choose $k_0$ large enough so that
$P_{k_0}\le(DM^2 k_0^{N})^{-1}$,
$(1+1/k_0)^N \le e^{\aca}$,
%
%
\begin{equation}
\label{eq541c} %
k_0 \ge D\bigl(c_2 \vee
(1/c_2)\bigr)s_0^{2h_T}h_T^{h_T +1}
\end{equation}
and $k_0\ge N_0$ all hold, where $M=e^{\aca h_T}$, $s_0$ is as in (\ref{eq532})
and $N_0$ is as in Lemma~\ref{lem391}. Setting $\ooo= Me^{\aca
k_0}P_{k_0}$ implies (\ref{eq541})
holds for $k = k_0,\ldots,k_0 + h_T$, which we will need in order to
begin our induction argument.

It follows from the definition of $\ooo$ and the first two conditions
on $k_0$ that
\[
\ooo DM e^{-\aca k} \le k^{-N} \qquad\mbox{for all } k\ge
k_0.
\]
Setting $q_k = \alpha D (\ooo M e^{-\aca k} )^{D-1}$, it follows from
this that
%
%
\begin{equation}
\label{eq541e} %
q_k \le k^{-(D-1)N} \qquad\mbox{for all }
k\ge k_0, %
\end{equation}
which we will use throughout the induction argument for (\ref{eq541}).
In order to follow the basic induction argument, the reader should keep
in mind (\ref{eq541}) and (\ref{eq541e}), without worrying much about
the other inequalities.\vadjust{\goodbreak}

%
%
%
%
In order to demonstrate the inequality in (\ref{eq541}) with $i=k$, we
proceed as outlined in the beginning of the subsection,
employing the partitions $\pi$ given in (\ref{eq3111}) and the
events $A_{\pi}$, on which a sequence of arrivals and departures
occurs in the first cycle that
induces the partition $\pi$. We define $\Pi_k$, as before, as the set
of all partitions with final element
$i_m=k$; here, the first element will be $i_0+1$, with $i_0=k_1$, where
$k_1=k-h_T$.
In the present setting, we will pay more attention than in Section~\ref{sec3} to the
length of each of the sets in a partition $\pi$, setting $h_{\ell} =
|G_{\ell}|$, for
$\ell=1,\ldots,m_{\pi}$, for the number of elements in the $\ell$th
set $G_{\ell}$ of the partition; one has
$h_T=\sum_{\ell=1}^{m_{\pi}}h_{\ell}$.

An important step in computing an upper bound for $P_k$ is
Proposition~\ref{prop551}, which is the
analog of Proposition~\ref{prop3121}. Rather than employing $L_{\ell}(s)$
as in the proof of
Proposition~\ref{prop3121} for the upper bound for a set in the
partition, we employ
%
%
\begin{equation}
\label{eq543} %
J_{k,h}(s) \stackrel{\mathrm{def}} {=}
e^{-q_k s}\sum_{i=h}^{\infty}(q_k
s)^i/i!. %
\end{equation}
%
%
%
%
%
The quantity $J_{k,h}(s)$ is the
probability of at least $h$ events occurring for a mean-$q_k(s)$
Poisson random variable, and
dominates the probability that, over the time interval $(0,s]$, at
least $h$ arrivals occur
for a cavity process $X^{\mathcal{H}}(\cdot)$ with $Z^{\mathcal
{H}}(0)\ge k_1$ and
$S^{\mathcal{H}}(0)\ge s$. This bound follows from the upper bound in
(\ref{eq234}), together with
the induction hypothesis (\ref{eq541}) and our definition of $M$.
%
%
%
\begin{prop}
\label{prop551}
Consider a family of JSQ networks with the same assumptions holding as
in Proposition~\ref{prop531}, except
that (\ref{eq532}) is not assumed. Suppose that the induction
assumption (\ref{eq541}) holds
for given $h_T$ and for $k_0 \ge N_0$, where $N_0$ is as in Lemma \ref
{lem391}. Then,
%
%
\begin{eqnarray}
\label{eq552} %
P_k &\le&3\sum
_{\pi\in\Pi_k}(3k)^{m_{\pi}-1} \prod_{\ell=1}^{m_{\pi}-1}
\biggl(\int_{0}^{\infty}J_{k,h_{\ell}}(s) F(ds)
\biggr)
\nonumber
\\[-8pt]
\\[-8pt]
\nonumber
&&\hspace*{28pt}{}\times \int_{0}^{\infty}(k+s)J_{k,h_{m_{{\pi}}}}(s)
F(ds). %
\end{eqnarray}
\end{prop}
\begin{pf}
One can reason similarly to the argument for (\ref{eq3126}), in the
proof of Proposition~\ref{prop3121},
by computing an upper bound on $E[V_k; A_{\pi}]$. Summation over $\pi
\in\Pi_k$ and application
of (\ref{eq233}) will then imply (\ref{eq552}). The assumption $k_0
\ge N_0$ is needed only to
absorb the term $N_0$ when applying Lemma~\ref{lem391}.

One argues inductively, repeating the argument for (\ref{eq3126}),
except for the substitution of
$J_{k,h_{\ell}}(s)$ for $L_{\ell}(s)$ and a minor change
involving the factors of $3k$. For each step with $\ell< m_{\pi}$,
one obtains an additional factor $i_{\ell-1}^*\int_{0}^{\infty
}J_{k,h_{\ell}}(s) F(ds)$ and, for
$\ell= m_{\pi}$, one obtains the factor
$3i_{m_{\pi} - 1}^* \int_{0}^{\infty}(k+s)J_{k,h_{m_{{\pi}}}}(s)
F(ds)$, where
$i_{\ell-1}^* = 3i_{\ell-1}$ for $\ell\ge2$ and $i_0^* = m_0$, with
$m_0$ being the mean return
time to $0$ for $X^{\mathcal{H}}(\cdot)$.
For $\ell< m_{\pi}$, the integral part of the factor
is obtained by employing the comparison given directly before
the statement of the proposition, comparing $J_{k,h_{\ell}}(s)$ with
the probability of at least $h$ arrivals
over a service time of at least $s$, and then by integrating against
$s$; for $\ell= m_{\pi}$, one also
employs (\ref{eq392}) to provide an upper bound on the expected
occupation time $V_k$.

For $\ell\ge2$, the factor $i_{\ell-1}^*$ is obtained by applying
(\ref{eq393}), with $s=0$, which gives an
upper bound on the expected number of service intervals occurring over
the remainder of the cycle,
after the service interval corresponding to the $(\ell-1)$st step ends.
For $\ell_0^*$, instead of the factor $3i_0$, one can employ
$m_0$, since this is the expected number of service
intervals over an entire cycle, and no conditioning is needed for this
first step. Since each of the
remaining factors is at most $3k$, the product of all of the factors is
at most $m_0(3k)^{m_{\pi}-1}$, and
since $P_k = (m_0)^{-1}E[V_k]$, the $m_0$ factors cancel, and one
obtains the $(3k)^{m_{\pi}-1}$ factor in
(\ref{eq552}). [The improved bound just obtained by removing a factor
of $3k$ will only be needed when
bounding the right-hand side of (\ref{eq552}) for the trivial
partition, in Proposition~\ref{prop561}.]
\end{pf}

In Propositions~\ref{prop561} and~\ref{prop571}, we provide upper
bounds for the summands on the right-hand side of (\ref{eq552}), which
we denote by $Q_k(\pi)$.
In Proposition~\ref{prop561}, we do this for the trivial partition consisting
of a single set, for which we write~$\pi_1$.
In Proposition~\ref{prop571}, we do this for each of the
other partitions. The sum over $\Pi-\{\pi_1\}$ of the bounds for
$Q_k(\pi)$ that are obtained
in Proposition~\ref{prop571} will be negligible
in comparison with the bound obtained for $Q_k(\pi_1)$
in Proposition~\ref{prop561}. This last bound will therefore dominate
the upper
bound for $P_k$ that will be obtained by inserting these bounds into
(\ref{eq552}) of the preceding proposition.

Both Propositions~\ref{prop561} and~\ref{prop571} employ the
elementary upper bounds for $J_{k,h}(s)$,
%
\begin{eqnarray}
\label{eq554} %
J_{k,h}(s) &\le&\bigl(4(q_k
s)^h /h!\bigr) \wedge1\qquad \mbox{for } s \le h/4q_k,
\nonumber
\\[-8pt]
\\[-8pt]
\nonumber
&\le&1\qquad \mbox{ for } s > h/4q_k, %
\end{eqnarray}
which one obtains by dominating the series in (\ref{eq543}) by the
geometric series $((q_k s)^{h}/h!)\sum_{i=0}^{\infty}(3/4)^i$, for $s
\le h/4q_k$.

%
\begin{prop}
\label{prop561}
%
Suppose that
%
%
\begin{equation}
\label{eq563} %
Q_k(\pi_1) = \int
_{0}^{\infty}(k+s)J_{k,h_T}(s) F(ds),
\end{equation}
where $F(\cdot)$ satisfies (\ref{eq532}) and $J_{k,h_T}(s)$ is chosen
as above, with
$h_T \ge6$, and suppose that
$k\ge k_0$, with (\ref{eq541c}) and~\ref{eq541e}) both holding. Then,
%
%
\begin{equation}
\label{eq562} %
Q_k(\pi_1) \le55D
c_2 (q_k/h_T)^{1/(D-1)}. %
\end{equation}
\end{prop}
\begin{pf}
%
%
%
%
%
Throughout the proof, we will abbreviate by setting $h_T=h$.
We begin the argument by decomposing the integral into the three parts,
$\int_0^{k}$, $\int_k^{h/27q_k}$ and $\int_{h/27q_k}^{\infty}$,
which we analyze separately.

Since $k+s\le2k$ for $s\in[0,k]$, it is easy to check that
%
%
\begin{equation}
\label{eq564} %
\int_{0}^{k}(k+s)J_{k,h}(s)
F(ds) \le8 k^{h+1}q_k^{h}/h!. %
\end{equation}

One has $k\ge s_0$ for $s_0$ in (\ref{eq532}). Applying (\ref{eq532})
and $k+s\le2s$,
and substituting $t=q_k s/h$, one sees that the second integral is
bounded above by
%
%
\begin{equation}
\label{eq565} \bigl(8D/(D-1)\bigr)c_2\int_{0}^{1/27}
q_k^{1/(D-1)}\frac{(t^{2/3}h)^h}{h!} t^{ ({h}/{3}-{D}/{(D-1)}
)} \,dt.
\end{equation}
Since $h\ge6$, one can check that $(t^{2/3}h)^h/h! \le3^{-h}$ and
$t^{ ({h}/{3}-{D}/{(D-1)} )} \le1$ for $t\le1/27$. Therefore,
(\ref{eq565}) is bounded above by
%
%
\begin{equation}
\label{eq566} %
(8/27) \bigl(D/(D-1)\bigr)c_2
3^{-h}q_k^{1/(D-1)} \le c_2
3^{-h}q_k^{1/(D-1)}. %
\end{equation}

Applying (\ref{eq532}), the third integral is at most
%
%
\begin{equation}
\label{eq567} 2\bigl(D/(D-1)\bigr)c_2\int_{h/27q_k}^{\infty}
s^{-D/(D-1)} \,ds \le54Dc_2(q_k/h)^{1/(D-1)}.
\end{equation}

On account of (\ref{eq541e}) and $q_k \le c_2$, the bound for the
third integral is clearly the dominant term.
Combining the bounds for the three integrals therefore implies that
\[
Q_k(\pi_1) \le55Dc_2 (q_k/h)^{1/(D-1)},
\]
which is the bound in (\ref{eq562}).
\end{pf}


%
\begin{prop}
\label{prop571}
%
%
Suppose that
%
%
\begin{eqnarray}
\label{eq573} %
Q_k(\pi) & = &(3k)^{m_{\pi}-1} \prod
_{\ell=1}^{m_{\pi}-1} \biggl(\int_{0}^{\infty}J_{k,h_{\ell}}(s)
F(ds) \biggr)
\nonumber
\\[-8pt]
\\[-8pt]
\nonumber
&&{}\times\int_{0}^{\infty}(k+s)J_{k,h_{m_{{\pi}}}}(s)
F(ds), %
\end{eqnarray}
where $F(\cdot)$ satisfies (\ref{eq532}) and $J_{k,h_{\ell}}(s)$ is
chosen as above, with
$h_T \ge5$, and suppose that
$k\ge k_0$, with (\ref{eq541c}) and (\ref{eq541e}) both holding. Then,
%
%
\begin{equation}
\label{eq572} %
Q_k(\pi) \le81 D^2(c_2
+1)^2 s_0^{2h_T} h_T^{h_T}
(3k)^{h_T}q_k^{D/(D-1)} %
\end{equation}
for each $\pi\in\Pi_k - \{\pi_1\}$.
\end{prop}

In order to demonstrate Proposition~\ref{prop571}, we will categorize
each partition in $\Pi_k - \{\pi_1\}$
as one of three types, based on the sizes and indices of its
constituent sets $G_{\ell}$, $\ell=1,\ldots,m_{\pi}$.
We will say $G_{\ell}$ is \emph{large} if $h_{\ell} \ge3$ and \emph
{small} if $h_{\ell} \le2$; we will also
distinguish between sets $G_{\ell}$ with $\ell< m_{\pi}$ and $\ell=
m_{\pi}$.
We will say that a partition $\pi$ is of type (I) if at least one
of its sets $G_{\ell}$, with $\ell< m_{\pi}$, is large; that it is
of type (II) if $G_{m_{\pi}}$ is large, but all
of the other sets are small; and that it is of type (III) if none of
its sets is large, but at least two sets
$G_{\ell_1}$ and $G_{\ell_2}$, with $\ell_1 < \ell_2 < m_{\pi}$
are small. It is easy to check that, for any
$h_T \ge5$, the three types of sets partition $\Pi_k - \{\pi_1\}$.

\begin{pf*}{Proof of Proposition~\ref{prop571}}
We will show separately that (\ref{eq572}) holds when~$\pi$ is a
member of any of the above three types of
partitions.
%
%
%
%
We will first bound the above integrals for the large and small sets
$G_{\ell}$, for both $\ell= m_{\pi}$ and
$\ell< m_{\pi}$, and will then apply these bounds to the three types
of partitions.
When convenient, we abbreviate by setting $h_{\ell}=h$.

Applying almost the same reasoning as in the proof of Proposition \ref
{prop561},
one obtains, for large $G_{m_{\pi}}$,
%
%
\begin{equation}
\label{eq581} %
\int_{0}^{\infty}(k+s)J_{k,h_{m_{{\pi}}}}(s)
F(ds) \le2D c_2 h_{m_{\pi}}^{h_{m_{\pi}}}
q_k^{1/(D-1)}. %
\end{equation}
One decomposes the integral into the parts $\int_{0}^{k}$, $\int
_{k}^{h/q_k}$ and $\int_{h/q_k}^{\infty}$.
A bound for the first integral is again given by the right-hand side
of (\ref{eq564}) and a bound for the third integral is
given by $2Dc_2(q_k/h)^{1/(D-1)}$. For the second integral, one obtains
the bound
$c_2h^h q_k^{1/(D-1)}$, after substituting $t=q_k s/h$ as before.
Instead of (\ref{eq565}),
one employs
%
%
\begin{equation}
\label{eq582} %
8\bigl(D/(D-1)\bigr)c_2 \int
_{0}^{1} q_k^{1/(D-1)}
\frac{h^h}{h!}t^{ (h-{D}/{(D-1)} )} \,dt %
\end{equation}
as an intermediate bound for the second integral, to which one applies
$t^{ (h-{D}/{(D-1)} )}\le1$;
the acquired factor $h^h$ will not cause difficulties in
the present context. For $k\ge k_0$, the bound in (\ref{eq581})
follows from the bounds
on the three integrals, on account of (\ref{eq541e}) and $q_k \le c_2$.

Similar reasoning can be applied for large $G_{\ell}$, with $\ell<
m_{\pi}$, to obtain the upper bound
%
%
\begin{equation}
\label{eq583} %
\int_{0}^{\infty}J_{k,h_{\ell}}(s)
F(ds) \le2 c_2 h_{\ell}^{h_{\ell}} q_k^{D/(D-1)}.
\end{equation}
One decomposes the integral into the parts $\int_{0}^{s_0}$, $\int
_{s_0}^{h/q_k}$ and $\int_{h/q_k}^{\infty}$.
The first integral is at most $s_0^h q_k^h \le s_0^h q_k^3$
and the third integral is at most $c_2 q_k^{D/(D-1)}$. For the
second integral, one obtains the upper bound $c_2 h^h q_k^{D/(D-1)}$,
after substituting $t=q_k s/h$.
Instead of (\ref{eq565}) or (\ref{eq582}), one employs
%
%
\begin{equation}
\label{eq584} %
4\bigl(D/(D-1)\bigr)c_2 \int
_{0}^{1} q_k^{D/(D-1)}
\frac{h^{h-1}}{h!}t^{ (h-({D}/{(D-1)})-1 )} \,dt\vadjust{\goodbreak} %
\end{equation}
as an intermediate bound for the second integral, to which one applies
$t^{ (h-({D}/{(D-1)})-1 )}\le1$.
Since $1/q_k \ge s_0^h$, the bound in (\ref{eq581}) follows from the
bounds on the three integrals.

For small $G_{\ell}$ with $\ell< m_{\pi}$, one obtains the upper bound
%
%
\begin{equation}
\label{eq585} %
\int_{0}^{\infty}J_{k,h_{\ell}}(s)
F(ds) \le9D(c_2 +1)s_0^{h_{\ell}}
q_k. %
\end{equation}
As in the previous case, one decomposes the integral into the parts
$\int_{0}^{s_0}$, $\int_{s_0}^{h/q_k}$ and $\int_{h/q_k}^{\infty}$.
The same estimates show that the first integral is at most $s_0^h q_k^h
\le s_0^h q_k$
and the third integral is at most $c_2 q_k^{D/(D-1)}$. For the second integral,
one obtains the upper bounds
%
%
\begin{equation}
\label{eq591} %
4\bigl(D/(D-1)\bigr)c_2 \int
_{s_0}^{h/q_k} q_k \frac{h^{h-1}}{h!}s^{-D/(D-1)}
\,ds \le8Dc_2 s_0 q_k, %
\end{equation}
with the inequality using $h\le2$. The bound in (\ref{eq585})
follows from the bounds on the three integrals.

For small $G_{m_{\pi}}$, the upper bound
%
%
\begin{equation}
\label{eq592} %
\int_{0}^{\infty}(k+s)J_{k,h_{m_{{\pi}}}}(s)
F(ds) \le k+1 \le2k %
\end{equation}
follows from $J_{k,h_{m_{{\pi}}}}(s) \le1$, since $F(\cdot)$ has
mean $1$.

We also note that, for $G_{\ell}$ with $\ell< m_{\pi}$,
%
%
\begin{equation}
\label{eq593} %
\int_{0}^{\infty}J_{k,h_{\ell}}(s)
F(ds) \le1 %
\end{equation}
trivially holds.

We now combine the upper bounds in (\ref{eq581}), (\ref{eq583}),
(\ref{eq585})
(\ref{eq592}) and (\ref{eq593}) to obtain upper bounds for the
right-hand side of (\ref{eq573}),
for large $k$.
When $\pi$ is a type~(I) partition, it follows from (\ref{eq583}),
(\ref{eq592}) and
(\ref{eq593}) that
%
%
\begin{equation}
\label{eq594} %
Q_k(\pi) \le2c_2
h_T^{h_T} (3k)^{m_{\pi}} q_k^{D/(D-1)};
\end{equation}
when $\pi$ is a type (II) partition, it follows from (\ref{eq581}),
(\ref{eq585}) and
(\ref{eq593}) that
%
%
\begin{equation}
\label{eq595} %
Q_k(\pi) \le18D^2(c_2+1)^2
s_0^{h_T} h_T^{h_T}
(3k)^{m_{\pi}-1} q_k^{D/(D-1)}; %
\end{equation}
and when $\pi$ is a type (III) partition, it follows from (\ref
{eq585}), (\ref{eq592}) and
(\ref{eq593}) that
%
%
\begin{equation}
\label{eq596} %
Q_k(\pi) \le81 D^2
(c_2+1)^2 s_0^{2h_T}
(3k)^{m_{\pi}} q_k^{2}. %
\end{equation}
The right-hand side of (\ref{eq572}) is greater than each of the
quantities in
(\ref{eq594})--(\ref{eq596}). Consequently, (\ref{eq572}) holds for all
$\pi\in\Pi_{k}-\{\pi_1\}$, as desired.
%
\end{pf*}

Employing Propositions~\ref{prop551},~\ref{prop561} and \ref
{prop571}, and the induction hypothesis
(\ref{eq541}), we now complete the proof of Proposition~\ref{prop531}.
\begin{pf*}{Proof of Proposition~\ref{prop531}}
We will demonstrate that the inequality in~(\ref{eq541}) holds for $i=k$,
provided it holds for $i=k_0,\ldots,k-1$, for $h_T$ and $\aca$ satisfying
(\ref{eq541a}) and (\ref{eq541b}), and for $k_0$ satisfying the inequalities
in (\ref{eq541c}) and on each side. By induction, it will follow that
%
%
\begin{equation}
\label{eq598} %
P_k\le\ooo e^{-\aca k} \qquad\mbox{for
all } k\ge k_0. %
\end{equation}

By Proposition~\ref{prop551},
%
%
\begin{equation}
\label{eq5101} %
P_k \le3\sum
_{\pi\in\Pi_k}Q_k(\pi) \le3Q_k(
\pi_1) + 3\cdot2^{h_T} \max_{\pi\in\Pi_k-\{\pi_1\}}Q_k(
\pi). %
\end{equation}
On account of (\ref{eq541c}) and (\ref{eq541e}), it follows from the
bounds in (\ref{eq562}) and (\ref{eq572}), for $Q_k(\pi_1)$ and for
$Q_k(\pi)$,
$\pi\in\Pi_k-\{\pi_1\}$, that the first term on the right-hand side
of (\ref{eq5101})
dominates the second term, and therefore
%
%
\begin{equation}
\label{eq5102} %
P_k \le220D c_2
(q_k/h_T)^{1/(D-1)}. %
\end{equation}
Substituting for $q_k$ and then for $M$, this is at most
%
%
\begin{equation}
\label{eq5103} %
\bigl(220 D^2 c_2
h_T^{-1/(D-1)} e^{\aca h_T} \bigr) \ooo
e^{-\aca k}. %
\end{equation}

Upon substitution of the value for $h_T$ in (\ref{eq541a}) and $\aca=
1/h_T$, the
quantity inside the parentheses in (\ref{eq5103}) is less than $1$.
Also, by
replacing the term $h_T^{-1/(D-1)}$ by $1$, it is easy to see that the
quantity inside the parentheses is again less than $1$, for
$\aca= \frac{1}{6}\log((220 D^2 c_2)^{-1} )$. So, in either case, the
inequality in (\ref{eq541}) holds for $i=k$. This implies (\ref{eq598}).

With a large enough choice of $\ooo$, (\ref{eq598}) extends to all
$k\ge0$. This
implies~(\ref{eq533}) of Proposition~\ref{prop531} with $r_D(c_2) =
\aca$, for this
choice of $\ooo$. Moreover, as $c_2\searrow0$, one has $\aca\nearrow
\infty$, and so (\ref{eq534})
also holds. This completes the proof of Proposition~\ref{prop531}.
\end{pf*}

%
%


\printaddresses

\end{document}